\documentclass[11pt,a4paper]{article}
\usepackage{amsmath,amssymb}
\usepackage{latexsym}
\usepackage{mathrsfs}
\usepackage{amsfonts}
\usepackage{hyperref}
\usepackage{amsthm}
\usepackage{appendix}
\usepackage[square,comma,sort&compress,numbers]{natbib}
\usepackage{color}
\usepackage{subcaption}
\usepackage{wrapfig}

\hfuzz=\maxdimen
\theoremstyle{definition}
\newtheorem{lemma}{Lemma}[section]
\newtheorem{definition}[lemma]{Definition}
\newtheorem{theorem}[lemma]{Theorem}
\newtheorem{corollary}[lemma]{Corollary}

\newtheorem{remark}{Remark}

\newtheorem{question}[lemma]{Question}
\newtheorem{conjecture}[lemma]{Conjecture}

\DeclareFixedFont{\Acknowledgment}{OT1}{cmr}{bx}{n}{14pt}
\begin{document}

\title{\bf  Parameterized combinatorial curvatures and parameterized combinatorial curvature flows
for discrete conformal structures on polyhedral surfaces}
\author{Xu Xu, Chao Zheng}
\maketitle

\begin{abstract}
Discrete conformal structure on polyhedral surfaces is a discrete analogue of the smooth conformal structure on surfaces
that assigns discrete metrics by scalar functions defined on vertices.
It unifies and generalizes tangential circle packing, Thurston's circle packing, inversive distance circle packing and
the vertex scaling described by Luo \cite{Luo 1} and others.
In this paper, we introduce combinatorial $\alpha$-curvature for discrete conformal structures on polyhedral surfaces,
which is a parameterized generalization of the classical combinatorial curvature.
Then we prove the local and global rigidity of combinatorial $\alpha$-curvature
with respect to discrete conformal structures on polyhedral surfaces,
which confirms parameterized Glickenstein rigidity conjecture in \cite{Glickenstein}.
To study the Yamabe problem for combinatorial $\alpha$-curvature,
we introduce combinatorial $\alpha$-Ricci flow for discrete conformal structures on polyhedral surfaces,
which is a generalization of Chow-Luo's combinatorial Ricci flow for Thurston's circle packings \cite{Chow-Luo} and
Luo's combinatorial Yamabe flow for vertex scaling \cite{Luo 1} on polyhedral surfaces.
To handle the potential singularities of the combinatorial $\alpha$-Ricci flow,
we extend the flow through the singularities by extending the inner angles in triangles by constants.
Under the existence of a discrete conformal structure with prescribed combinatorial curvature,
the solution of extended combinatorial $\alpha$-Ricci flow is proved to exist for all time and
converge exponentially fast for any initial value.
This confirms a parameterized generalization of another conjecture of Glickenstein in \cite{Glickenstein}  on the convergence of combinatorial Ricci flow, gives an almost equivalent characterization of  the solvability of  Yamabe problem for combinatorial $\alpha$-curvature
in terms of combinatorial $\alpha$-Ricci flow
and provides an effective algorithm for finding discrete conformal structures with prescribed combinatorial $\alpha$-curvatures.
\end{abstract}


\textbf{Keywords}: Combinatorial Ricci flow; Rigidity; Discrete conformal structure; Combinatorial curvature

\section{Introduction}
This is a continuation of \cite{Xu 21} studying discrete conformal structures on polyhedral surfaces.
Discrete conformal structure on polyhedral surfaces is a discrete analogue of the smooth conformal structure on Riemannian surfaces,
which defines discrete metrics by scalar functions defined on the vertices.
There are mainly four special types of discrete conformal structures on polyhedral surfaces that have been extensively studied in history,
including the tangential circle packing, Thurston's circle packing, inversive distance circle packing and the vertex scaling described by Luo \cite{Luo 1} and others.
These special types of discrete conformal structures are introduced and studied individually in the past.
The generic notion of discrete conformal structure on polyhedral surfaces was
introduced recently independently by Glickenstein \cite{Glickenstein2011}  and Glickenstein-Thomas \cite{Glickenstein-Thomas} from Riemannian geometry perspective
and by Zhang-Guo-Zeng-Luo-Yau-Gu \cite{ZGZLYG} using 3-dimensional hyperbolic geometry,
which unifies and generalizes the existing special types of discrete conformal structures on polyhedral surfaces.
In \cite{Xu 21}, the first author studied the classical combinatorial curvature for discrete conformal structures on polyhedral surfaces.
In this paper, we introduce a new combinatorial curvature (combinatorial $\alpha$-curvature) for discrete conformal structures on polyhedral surfaces, which is a parameterized generalization of the classical combinatorial curvature.
Then we prove its rigidity and further study the corresponding Yamabe problem for the new combinatorial curvature
using combinatorial Ricci flow and combinatorial Calabi flow.

Suppose $(M, \mathcal{T})$ is a triangulated connected closed surface with a triangulation $\mathcal{T}={(V,E,F)}$, where $V,E,F$ represent the sets of vertices, edges and faces respectively and $V$ is a finite subset of $M$ with $|V|=N$.
Denote a vertex, an edge and a face in the triangulation $\mathcal{T}$ by $i, \{ij\}, \{ijk\}$ respectively, where ${i,j,k}$ are natural numbers.
If a map $l: E\rightarrow (0, +\infty)$ assigns a length to every edge in such a way
that every triangle $\{ijk\}\in F$ with edge lengths $l_{ij}, l_{ik}, l_{jk}$ is embedded in 2-dimensional Euclidean space $\mathbb{E}^2$ (2-dimensional hyperbolic space $\mathbb{H}^2$ or 2-dimensional spherical space $\mathbb{S}^2$ respectively),
then $(M, \mathcal{T}, l)$ is called
as a triangulated polyhedral surface with Euclidean (hyperbolic or spherical respectively) background geometry and $l: E\rightarrow (0, +\infty)$ is called a Euclidean (hyperbolic or spherical respectively) polyhedral metric.
One can also take the triangulated polyhedral surface $(M, \mathcal{T}, l)$ being obtained
by gluing triangles in $\mathbb{E}^2$ ($\mathbb{H}^2$ or $\mathbb{S}^2$ respectively) isometrically along the edges in pair.
For a triangulated polyhedral surface $(M, \mathcal{T}, l)$, the classical combinatorial curvature $K: V\rightarrow (-\infty, 2\pi)$ is used to describe the conic singularities of polyhedral metrics at the vertices, which is defined to be
\begin{equation}\label{classical combinatorial curvature}
K_i=2\pi-\sum_{\{ijk\}\in F}\theta_i^{jk}
\end{equation}
with summation taken over all the triangles with $i$ as a vertex and $\theta_i^{jk}$ being the inner angle of the triangle $\{ijk\}\in F$ at the vertex $i$.
The classical combinatorial curvature satisfies the following discrete Gauss-Bonnet formula (\cite{Chow-Luo}, Proposition 3.1)
\begin{equation*}
\sum_{i\in V}K_i=2\pi\chi(M)-\lambda Area(M),
\end{equation*}
where $\lambda=-1, 0, +1$ for hyperbolic, Euclidean and spherical background geometry respectively and $Area(M)$ is the area of the surface $M$.

\begin{definition}(\cite{Glickenstein-Thomas,ZGZLYG})\label{edge length definition}
Suppose $(M,\mathcal{T})$ is a triangulated connected closed surface with two weights $\varepsilon: V\rightarrow \{-1,0,1\}$ and $\eta: E\rightarrow \mathbb{R}$ satisfying $\eta_{ij}=\eta_{ji}$.
A discrete conformal structure on $(M, \mathcal{T}, \varepsilon, \eta)$ is a map $f: V\rightarrow \mathbb{R}$
determining a polyhedral metric $l:E\rightarrow \mathbb{R}$
by assigning the edge length $l_{ij}$ for $\{ij\}\in E$ as
\begin{equation}\label{Euclidean edge length}
l_{ij}=\sqrt{\varepsilon_ie^{2f_i}+\varepsilon_je^{2f_j}+2\eta_{ij}e^{f_i+f_j}}
\end{equation}
in the Euclidean background geometry,
\begin{equation}\label{hyperbolic edge length}
l_{ij}=\cosh^{-1}\left(\sqrt{(1+\varepsilon_ie^{2f_i})(1+\varepsilon_je^{2f_j})}+\eta_{ij}e^{f_i+f_j}\right)
\end{equation}
in the hyperbolic background geometry and
\begin{equation*}
l_{ij}=\cos^{-1}\left(\sqrt{(1-\varepsilon_ie^{2f_i})(1-\varepsilon_je^{2f_j})}-\eta_{ij}e^{f_i+f_j}\right)
\end{equation*}
in the spherical background geometry.
The weight $\varepsilon: V\rightarrow \{-1,0, 1\}$ is called the scheme coefficient and $\eta: E\rightarrow \mathbb{R}$ is called the discrete conformal structure coefficient.
\end{definition}
Two discrete conformal structures defined on the same weighted triangulated surface $(M,\mathcal{T},\varepsilon,\eta)$ with the same background geometry is said to be conformally equivalent. In this paper, we focus on the cases of Euclidean and hyperbolic background geometry.

\begin{remark}
The relationships of the discrete conformal structure in Definition \ref{edge length definition} and the existing main special types of discrete conformal structures are summarized in Table \ref{Relationships of DCS}.
\begin{table}[ht]
\centering
\begin{tabular}{|l l l l|}
\hline
  Scheme & $\varepsilon_i$ & $\varepsilon_j$ & $\eta_{ij}$ \\ \hline
  Tangential circle packing & $+1$ & $+1$ & $+1$ \\
  Thurston's circle packing & $+1$ & $+1$ & $(-1,1]$ \\
  Inversive distance circle packing & $+1$ & $+1$ & $(-1, +\infty)$ \\
  Vertex scaling & $0$ & $0$ & $(0, +\infty)$ \\
  Discrete conformal structure & $\{+1, 0, -1\}$ & $\{+1, 0, -1\}$ & $(-1, +\infty)$ \\
  \hline
\end{tabular}
\caption{Relationships of discrete conformal structures} 
\label{Relationships of DCS} 
\end{table}
By Table \ref{Relationships of DCS}, the discrete conformal structure in Definition \ref{edge length definition} not only contains vertex scaling and different types of  circle packings as special cases, but also contains other types of discrete conformal structures, for example the mixed type with $\varepsilon_i=0$ for vertices $i\in V_0\neq \emptyset$ and $\varepsilon_j=1$ for the other vertices $j\in V\setminus V_0\neq \emptyset$.
Please refer to \cite{Glickenstein-Thomas, ZGZLYG} for more information on this.
\end{remark}

The classical combinatorial curvature $K$ on polyhedral surfaces is proved to converge to the smooth Gaussian curvature on surfaces
in the sense of measure by Cheeger-M\"{u}ller-Schrader \cite{CMS}. However, the convergence
is not a purely local phenomenon.
As mentioned in \cite{CMS,Ge-Xu 21}, the classical combinatorial curvature $K$ defined by (\ref{classical combinatorial curvature}) is scaling invariant in the Euclidean background geometry and
does not approximate the smooth Gaussian curvature pointwisely on smooth surfaces as the triangulation of the surface becomes finer and finer.
The key issue is how to decompose such a measure over the triangulated surfaces.
Motivated by this idea and the original definition of smooth Gaussian curvature via Gauss map on surfaces,
Ge and the first author \cite{Ge-Xu 21} introduced the combinatorial $\alpha$-curvature for Thurston's Euclidean circle packing metrics
on triangulated surfaces.
After that, there are lots of works on new combinatorial curvatures on polyhedral surfaces and triangulated $3$-manifolds. Please refer to \cite{Dai-Ge,Ge-Jiang III,Ge-Xu alpha-curvatures,Ge-Xu 14,Ge-Xu 21,Xu 18,Xu 20,Ge-Xu 18,Ge-Xu 17, Kourimska thesis,Kourimska paper, Xu I,Xu Z,Xu Z 2} and other works.
In this paper, we introduce the following combinatorial $\alpha$-curvature for discrete conformal structures on a weighted triangulated surface $(M,\mathcal{T},\varepsilon,\eta)$, which is a parameterized generalization of the classical combinatorial curvature $K$.
\begin{definition}\label{combinatorial alpha_curvature definition}
Suppose $(M,\mathcal{T},\varepsilon,\eta)$ is a weighted triangulated
connected closed surface with weights $\varepsilon: V\rightarrow \{0,1\}$ and $\eta: E\rightarrow \mathbb{R}$.
The combinatorial $\alpha$-curvature at the vertex $i\in V$ is defined to be
\begin{equation*}
R_{\alpha, i}=\frac{K_i}{e^{\alpha u_i}},
\end{equation*}
where $K_i$ is the classical combinatorial curvature at $i\in V$ defined by (\ref{classical combinatorial curvature}),
$u_i=f_i$ for $i\in V$ in the Euclidean background geometry and
\begin{eqnarray}\label{hyperblic DCS in u}
u_i=
\begin{cases}
f_i,  &{\varepsilon_i=0},\\
\frac{1}{2}\ln\left|\frac{\sqrt{1+e^{2f_i}}-1}{\sqrt{1+e^{2f_i}}+1}\right|, &{\varepsilon_i=1}
\end{cases}
\end{eqnarray}
for $i\in V$ in the hyperbolic background geometry.
\end{definition}
We also call the function $u$ in Definition \ref{combinatorial alpha_curvature definition}
as a discrete conformal structure on $(M,\mathcal{T},\varepsilon,\eta)$,
if it causes no confusion in the context.

\begin{remark}
The combinatorial $\alpha$-curvature $R_{\alpha}$  in Definition \ref{combinatorial alpha_curvature definition}
is motivated by the following facts.
The first is that the smooth Gaussian curvature on surfaces is defined to be the limit of area distortion of the Gauss map.
In the Euclidean background geometry,
if we set $e^{f_i}=r_i$, which is taken to be the radius of the disk attached to the vertex $i\in V$, then
$e^{2u_i}=r_i^2$ is the area of the disk with radius $r_i$ up to a constant.
Similarly, in the hyperbolic background geometry, set $e^{f_i}=\sinh r_i$, then
$e^{2u_i}$ is a 1-st order approximation of the area of hyperbolic disk with radius $r_i$ up to a positive constant.
In this sense, the combinatorial $2$-curvature could be taken as an approximation of the smooth Gaussian curvature and
the combinatorial $\alpha$-curvature $R_{\alpha}$ is a parameterized generalization of it.
The second is that the matrix $(\frac{\partial \theta}{\partial u})$ is symmetric by Lemma \ref{Euclidean extended lemma} and
Lemma \ref{hyperbolic extended lemma}. This provides lots of convenience in the arguments.
The third is that this form of combinatorial $\alpha$-curvature unifies and generalizes
the previously defined combinatorial curvatures on polyhedral surfaces
and preserves similar scaling property of the smooth Gaussian curvature on surfaces.
If $\alpha=0$, the combinatorial curvature  $R_{\alpha}$  in Definition \ref{combinatorial alpha_curvature definition} is reduced to the classical combinatorial curvature $K$.
If $\varepsilon\equiv 1$, the combinatorial curvature $R_\alpha$ in Definition \ref{combinatorial alpha_curvature definition}
is reduced to the parameterized combinatorial curvature for circle packings studied in
\cite{Ge-Xu alpha-curvatures, Ge-Jiang III,Ge-Xu 14,Ge-Xu 17,Ge-Xu 18,Ge-Xu 21,Xu 18}.
If $\varepsilon\equiv 0$, the combinatorial curvature $R_\alpha$ in Definition \ref{combinatorial alpha_curvature definition}
is reduced to the parameterized combinatorial curvature for vertex scaling studied in \cite{Xu I,Xu Z,Xu Z 2}.
The parameterized combinatorial curvature $R_\alpha$ in Definition \ref{combinatorial alpha_curvature definition}
further covers the case of mixed type discrete conformal structures, for example, $\varepsilon_i=0$ for vertices $i\in V_0\neq \emptyset$ and $\varepsilon_j=1$ for the other vertices $j\in V\setminus V_0\neq \emptyset$.
\end{remark}

\begin{remark}\label{remark on u}
In the Euclidean background geometry, if we take $g_i=e^{u_i}$ as a discrete analogue of the smooth Riemannian metric,
then for any constant $\lambda>0$, we have
$R_{\alpha, i}(\lambda g_1,\cdots, \lambda g_N)=\lambda^{-\alpha}R_{\alpha, i}( g_1,\cdots, g_N)$.
In the special case of $\alpha=1$, $R_{1, i}(\lambda g_1,\cdots, \lambda g_N)=\lambda^{-1}R_{1, i}(g_1,\cdots, g_N)$, which is parallelling to the transformation of smooth Gaussian curvature $K_{\lambda g}=\lambda^{-1} K_g$ with $g$ being the Riemannian metric.
Note that, in the hyperbolic background geometry, $u_i\in \mathbb{R}$ for vertex $i$ with $\varepsilon_i=0$ and $u_i\in \mathbb{R}_{<0}$ for vertex $i$ with $\varepsilon_i=1$, thus $u=(u_1,..,u_N)\in \mathbb{R}^{N_0}\times \mathbb{R}^{N_1}_{<0}$, where $N_0$ is the number of  vertices $i\in V$ with $\varepsilon_i=0$ and $N_1=N-N_0$.
\end{remark}

A basic problem in discrete conformal geometry is to understand the relationships between the discrete conformal structure
and its combinatorial curvature.
For the combinatorial $\alpha$-curvature $R_\alpha$, we have the following rigidity with respect to discrete conformal structures on polyhedral surfaces, which confirms a parameterized generalization of Glickenstein's rigidity conjecture in \cite{Glickenstein}.
\begin{theorem}\label{global rigidity theorem}
Suppose $(M,\mathcal{T},\varepsilon,\eta)$ is a weighted triangulated connected closed surface with the weights $\varepsilon: V\rightarrow \{0,1\}$ and $\eta: E\rightarrow \mathbb{R}$ satisfying
\begin{equation}\label{structure condition 1}
\varepsilon_s\varepsilon_t+\eta_{st}>0,\ \forall\{st\}\in E
\end{equation}
and
\begin{equation}\label{structure condition 2}
\varepsilon_q\eta_{st}+\eta_{qs}\eta_{qt}\geq0,\ \{q,s,t\}=\{i,j,k\}
\end{equation}
for any triangle $\{ijk\}\in F$. $\alpha\in \mathbb{R}$ is a constant and $\overline{R}: V\rightarrow\mathbb{R}$ is a given function.
\begin{description}
  \item[(a)] In the case of Euclidean background geometry, if $\alpha\overline{R}\equiv 0$, there exists at most one discrete conformal structure $f: V\rightarrow \mathbb{R}$ with combinatorial $\alpha$-curvature $\overline{R}$ up to a vector $c(1,1,...,1),\ c\in \mathbb{R}$; if $\alpha\overline{R}\leq 0$ and $\alpha\overline{R}\not\equiv 0$, there exists at most one discrete conformal structure $f: V\rightarrow \mathbb{R}$ with combinatorial $\alpha$-curvature $\overline{R}$.
  \item[(b)] In the case of hyperbolic background geometry, if $\alpha\overline{R}\leq 0$, there exists at most one discrete conformal structure $f: V\rightarrow \mathbb{R}$ with combinatorial $\alpha$-curvature $\overline{R}$.
\end{description}
\end{theorem}
\begin{remark}
In the special case of $\alpha=0$, Theorem \ref{global rigidity theorem} was proved by the first author in \cite{Xu 21}.
If $\varepsilon_i=1$ for all $i\in V$, Theorem \ref{global rigidity theorem} is reduced to the rigidity of circle packings on surfaces obtained in \cite{Chow-Luo,Ge-Jiang III, Ge-Xu 18,Ge-Xu 21,Ge-Xu alpha-curvatures, Guo, Luo3, Xu 18} and others.
If $\varepsilon_i=0$ for all $i\in V$, Theorem \ref{global rigidity theorem} is reduced to the rigidity of vertex scaling on polyhedral surfaces obtained in \cite{BPS, Luo 1, Xu I,Xu Z} and others.
Theorem \ref{global rigidity theorem} unifies these results and further covers the case of mixed type that $\varepsilon_i=0$ for some vertices $i\in V_0\neq\emptyset$ and $\varepsilon_i=1$ for the other vertices $j\in V\setminus V_0\neq\emptyset$.
\end{remark}

It is interesting to consider the following combinatorial Yamabe problem for combinatorial $\alpha$-curvature $R_\alpha$.

\noindent $\mathbf{Combinatorial}$ $\mathbf{Yamabe\  problem \ }$ Does there exist a discrete conformal structure with constant combinatorial $\alpha$-curvature or prescribed combinatorial $\alpha$-curvature on $(M,\mathcal{T},\varepsilon,\eta)$? Furthermore, how to find it?

To study the combinatorial Yamabe problem for combinatorial $\alpha$-curvature $R_\alpha$,
we introduce the following combinatorial $\alpha$-Ricci flow for discrete conformal structures on polyhedral surfaces,
which is a generalization of Chow-Luo's combinatorial Ricci flow for Thurston's circle packings in \cite{Chow-Luo},
Luo's combinatorial Yamabe flow for vertex scaling in \cite{Luo 1} and
Zhang-Guo-Zeng-Luo-Yau-Gu's combinatorial Ricci flow for discrete conformal structures in \cite{ZGZLYG}.

\begin{definition}\label{Ricci flow definition}
Suppose $(M,\mathcal{T},\varepsilon,\eta)$ is a weighted triangulated connected closed surface with weights $\varepsilon: V\rightarrow \{0,1\}$ and $\eta: E\rightarrow \mathbb{R}$. $\alpha\in \mathbb{R}$ is a constant and $u: V\rightarrow \mathbb{R}$ is a discrete conformal structure defined by $u_i=f_i$ in the Euclidean background geometry and by (\ref{hyperblic DCS in u}) in the hyperbolic background geometry. The combinatorial $\alpha$-Ricci flow for discrete conformal structures on polyhedral surfaces is defined to be
\begin{equation}\label{combinatorial alpha Ricci flow}
\frac{du_i}{dt}=-R_{\alpha,i}
\end{equation}
for Euclidean and hyperbolic background geometry.
\end{definition}

\begin{remark}
If $\alpha=0$, the combinatorial $\alpha$-Ricci flow in Definition \ref{Ricci flow definition} is Zhang-Guo-Zeng-Luo-Yau-Gu's combinatorial Ricci flow in \cite{ZGZLYG}, which unifies and generalizes Chow-Luo's combinatorial Ricci flow for circle packing in \cite{Chow-Luo,Ge-Jiang I,Ge-Jiang II} and Luo's combinatorial Yamabe flow for vertex scaling in \cite{Luo 1}.
If $\varepsilon\equiv1$, the combinatorial $\alpha$-Ricci flow in Definition \ref{Ricci flow definition} is reduced to the
combinatorial $\alpha$-Ricci flow for circle packings studied in \cite{Ge-Jiang III,Ge-Xu 14,Ge-Xu 17,Ge-Xu 18,Ge-Xu 21,Ge-Xu alpha-curvatures}.
If $\varepsilon\equiv0$, the combinatorial $\alpha$-Ricci flow in Definition \ref{Ricci flow definition} is reduced to the
combinatorial $\alpha$-Yamabe flow for vertex scaling studied in \cite{Xu I,Xu Z}. The combinatorial $\alpha$-Ricci flow in Definition \ref{Ricci flow definition} further covers the case of mixed type discrete conformal structures.
\end{remark}

In the Euclidean background geometry,
the normalized combinatorial $\alpha$-Ricci flow is defined to be
\begin{equation}\label{normalized combinatorial alpha Ricci flow}
\frac{du_i}{dt}=R_{\alpha,av}-R_{\alpha,i},
\end{equation}
where $R_{\alpha,av}=\frac{2\pi \chi(M)}{\sum^N_{i=1}e^{\alpha u_i}}$ is the average combinatorial $\alpha$-curvature.
In the Euclidean and hyperbolic background geometry, we usually generalize the combinatorial $\alpha$-Ricci flow (\ref{combinatorial alpha Ricci flow})
to the following modified combinatorial $\alpha$-Ricci flow
\begin{equation}\label{modified combinatorial alpha Ricci flow}
\frac{du_i}{dt}=\overline{R}_i-R_{\alpha,i},
\end{equation}
where $\overline{R}$ is a given function defined on the vertices.

Along the combinatorial $\alpha$-Ricci flows (\ref{normalized combinatorial alpha Ricci flow}) and (\ref{modified combinatorial alpha Ricci flow}), singularities may develop, which correspond to that some triangle degenerates
or the discrete conformal structure tends to infinity.
If some triangle degenerates along the combinatorial curvature flow, we call the combinatorial curvature flow develops a \textit{removable singularity}.
If the discrete conformal structure tends to infinity, we call the the combinatorial curvature flow develops an \textit{essential singularity}.
To handle the potential removable singularities along the combinatorial $\alpha$-Ricci flows, we extend the combinatorial $\alpha$-curvature and then extend the combinatorial $\alpha$-Ricci flows through the removable singularities.
We have the following result on the longtime existence and convergence for the solution of the combinatorial $\alpha$-Ricci flow, which confirms a parameterized generalization of another conjecture of Glickenstein in \cite{Glickenstein} on the convergence of combinatorial Ricci flow, gives an almost equivalent characterization of  the solvability of  combinatorial Yamabe problem for combinatorial $\alpha$-curvature in terms of extended combinatorial $\alpha$-Ricci flow and provides effective algorithms for finding discrete conformal structures with constant combinatorial $\alpha$-curvature or prescribed combinatorial $\alpha$-curvatures.

\begin{theorem}\label{main theorem 1}
Suppose $(M,\mathcal{T},\varepsilon,\eta)$ is a weighted triangulated connected closed surface with the weights $\varepsilon: V\rightarrow \{0,1\}$ and $\eta: E\rightarrow \mathbb{R}$ satisfying the structure conditions (\ref{structure condition 1}) and (\ref{structure condition 2}). $\alpha\in \mathbb{R}$ is a constant and $\overline{R}: V\rightarrow \mathbb{R}$ is a given function.
\begin{description}
  \item[(a)] The solution of the combinatorial $\alpha$-Ricci flows (\ref{normalized combinatorial alpha Ricci flow})
  and  (\ref{modified combinatorial alpha Ricci flow}) could be extended by extending the inner angles in triangles by constants.
  Furthermore,  if $\alpha\chi(M)\leq0$,  the extended solution of the normalized combinatorial $\alpha$-Ricci flow (\ref{normalized combinatorial alpha Ricci flow}) is unique for any initial Euclidean discrete conformal structure;
  If $\alpha\overline{R}\leq 0$, the extended solution of modified combinatorial $\alpha$-Ricci flow (\ref{modified combinatorial alpha Ricci flow}) is unique for any initial Euclidean and hyperbolic discrete conformal structure.
  \item[(b)] Suppose there exists a Euclidean discrete conformal structure with constant combinatorial $\alpha$-curvature and $\alpha\chi(M)\leq0$. Then the normalized combinatorial $\alpha$-Ricci flow (\ref{normalized combinatorial alpha Ricci flow}) develops no essential singularities.
      If the solution of (\ref{normalized combinatorial alpha Ricci flow}) develops no removable singularities in finite time, then the solution of (\ref{normalized combinatorial alpha Ricci flow}) exists for all time, converges exponentially fast for any initial Euclidean discrete conformal structure and does not develop removable singularities at time infinity.
      Furthermore, the extended solution of the normalized combinatorial $\alpha$-Ricci flow (\ref{normalized combinatorial alpha Ricci flow}) in the Euclidean background geometry exists for all time and converges exponentially fast for any initial Euclidean discrete conformal structure.
  \item[(c)] Suppose there exists a hyperbolic discrete conformal structure with combinatorial $\alpha$-curvature $\overline{R}$ satisfying one of the following three conditions
\begin{description}
  \item[(1)] $\alpha>0$ and $\overline{R}_i\leq0$ for all $i\in V$,
  \item[(2)] $\alpha<0$ and $\overline{R}_i\in [0,2\pi)$ for all $i\in V$,
  \item[(3)] (\cite{Xu 21}) $\alpha=0$, $\overline{R}_i\in (-\infty, 2\pi)$ for all $i\in V$ and $\sum^N_{i=1} \overline{R}_{i}>2\pi \chi(M)$.
\end{description}
Then the modified combinatorial $\alpha$-Ricci flow (\ref{modified combinatorial alpha Ricci flow})
in the hyperbolic background geometry develops no essential singularities.
If the solution of (\ref{modified combinatorial alpha Ricci flow}) develops no removable singularities in finite time, then the solution of (\ref{modified combinatorial alpha Ricci flow}) exists for all time, converges exponentially fast for any initial hyperbolic discrete conformal structure and does not develop removable singularity at time infinity.
Furthermore, the extended solution of  modified combinatorial $\alpha$-Ricci flow (\ref{modified combinatorial alpha Ricci flow}) in the hyperbolic background geometry exists for all time and converges exponentially fast for any initial hyperbolic discrete conformal structure.
\end{description}
\end{theorem}
\begin{remark}
If $\alpha=0$, Theorem \ref{main theorem 1} can be extended to a form that holds for any reasonable prescribed combinatorial curvature with Euclidean and hyperbolic background geometry, which was proved by the first author in \cite{Xu 21}.
The result in Theorem \ref{main theorem 1} (b) can also generalized to the case of prescribed combinatorial $\alpha$-curvature in the Euclidean background geometry. Please refer to Theorem \ref{new theorem} for details.
The idea of extension to handle the singularities of the combinatorial $\alpha$-Ricci flow comes from Bobenko-Pinkall-Springborn \cite{BPS} and Luo \cite{Luo3}.
There is another approach to extend the combinatorial $\alpha$-Ricci flow for vertex scaling on polyhedral surfaces introduced in \cite{Gu1,Gu2}, which is called surgery by flipping under the Delaunay condition.
In the approach in \cite{Gu1,Gu2}, the condition on the existence of discrete conformal structure with  prescribed combinatorial curvature is removed. The approach in \cite{Gu1,Gu2} was recently used to study combinatorial curvature flows for vertex scaling on polyhedral surfaces in \cite{Zhu Xu,Xu I,Xu Z}.
\end{remark}

Combinatorial Calabi flow is another effective combinatorial curvature flow to study combinatorial Yamabe problem, which was first introduced by Ge \cite{Ge1} (see also \cite{Ge2}) for Thurston's Euclidean circle packing and then further studied in \cite{Ge3,Ge-Xu 16, Ge-Xu 21,Ge-Xu 14,GXZ,Xu I,Xu JFA 21,Xu 21,Xu Z, WX,Zhu Xu} and others.
We introduce the following combinatorial $\alpha$-Calabi flow for discrete conformal structures on polyhedral surfaces.

\begin{definition}\label{Calabi flow definition}
Suppose $(M,\mathcal{T},\varepsilon,\eta)$ is a weighted triangulated  connected closed surface with the weights $\varepsilon: V\rightarrow \{0,1\}$ and $\eta: E\rightarrow \mathbb{R}$. $\alpha\in \mathbb{R}$ is a constant and  $u: V\rightarrow \mathbb{R}$ is a discrete conformal structure defined by $u_i=f_i$ in Euclidean background geometry and by (\ref{hyperblic DCS in u}) in the hyperbolic background geometry. The combinatorial $\alpha$-Calabi flow for discrete conformal structures on polyhedral surfaces is defined to be
\begin{equation}\label{combinatorial alpha Calabi flow}
\frac{du_i}{dt}=\Delta_\alpha R_{\alpha,i}
\end{equation}
for Euclidean and hyperbolic background geometry,
where the discrete $\alpha$-Laplace operator $\Delta_\alpha$ is defined by
\begin{equation*}
\Delta_\alpha g_i=-\frac{1}{e^{\alpha u_i}}\sum_{j\in V}\frac{\partial K_i}{\partial u_j}g_j
\end{equation*}
for $g: V\rightarrow \mathbb{R}$.
\end{definition}

\begin{remark}
If $\alpha=0$, the combinatorial $\alpha$-Calabi flow in Definition \ref{Calabi flow definition} is the combinatorial Calabi flow introduced in \cite{Xu 21}, which unifies combinatorial Calabi flow for circle packings studied in \cite{Ge-Xu 16,Ge1,Ge2, Ge3} and combinatorial Calabi flow for vertex scaling studied in \cite{Ge1,Zhu Xu} and generalizes them to a very general context.
If $\varepsilon\equiv1$, the combinatorial $\alpha$-Calabi flow in Definition \ref{Calabi flow definition} is reduced to the
combinatorial $\alpha$-Calabi flow for circle packings studied in \cite{Ge-Xu 14,Ge-Xu 21}.
If $\varepsilon\equiv0$, the combinatorial $\alpha$-Calabi flow in Definition \ref{Calabi flow definition} is reduced to the
combinatorial $\alpha$-Calabi flow for vertex scaling studied in \cite{Xu I,Xu Z}. The combinatorial $\alpha$-Calabi flow in Definition \ref{Calabi flow definition} further covers the case of mixed type discrete conformal structures.
\end{remark}

We have the following result on the longtime existence and convergence for the solution of combinatorial $\alpha$-Calabi flow (\ref{combinatorial alpha Calabi flow}).
\begin{theorem}\label{main theorem 2}
Suppose $(M,\mathcal{T},\varepsilon,\eta)$ is a weighted triangulated  connected closed surface with the weights $\varepsilon: V\rightarrow \{0,1\}$ and $\eta: E\rightarrow \mathbb{R}$ satisfying the structure conditions (\ref{structure condition 1}) and (\ref{structure condition 2}). $\alpha\in \mathbb{R}$ is a constant and $\overline{R}: V\rightarrow \mathbb{R}$ is a given function.
\begin{description}
  \item[(a)] In the case of Euclidean background geometry, if the solution of the combinatorial $\alpha$-Calabi flow (\ref{combinatorial alpha Calabi flow}) converges to a nondegenerate discrete conformal structure, then there exists a discrete conformal structure
       on $(M,\mathcal{T},\varepsilon,\eta)$ with constant combinatorial $\alpha$-curvature. Furthermore, suppose that there exists a Euclidean discrete conformal structure $\overline{u}$ on $(M,\mathcal{T},\varepsilon,\eta)$ with constant combinatorial $\alpha$-curvature and $\alpha\chi(M)\leq 0$, then there exists a constant $\delta >0$ such that if
       the initial value $u(0)$ satisfies
      $||R_\alpha(u(0))-R_\alpha(\overline{u})||<\delta$ and $\sum^N_{i=1}e^{\alpha u(0)}=\sum^N_{i=1}e^{\alpha\overline{u}_i}$ in the case of $\alpha\neq 0$ or $||u(0)-\overline{u}||<\delta$ and $\sum^N_{i=1}u(0)=\sum^N_{i=1}\overline{u}_i$ in the case of $\alpha=0$,
       the solution of Euclidean combinatorial $\alpha$-Calabi flow (\ref{combinatorial alpha Calabi flow}) exists for all time and converges exponentially fast to $\overline{u}$.
  \item[(b)] In the case of hyperbolic background geometry, if the solution of the combinatorial $\alpha$-Calabi flow (\ref{combinatorial alpha Calabi flow}) converges to a nondegenerate discrete conformal structure, then there exists a discrete conformal structure on $(M,\mathcal{T},\varepsilon,\eta)$ with zero combinatorial $\alpha$-curvature. Furthermore, suppose that there exists a hyperbolic discrete conformal structure $\overline{u}$ on $(M,\mathcal{T},\varepsilon,\eta)$ with zero combinatorial $\alpha$-curvature, then there exists a constant $\delta >0$ such that if $||R_\alpha(u(0))-R_\alpha(\overline{u})||<\delta$, the solution of hyperbolic combinatorial $\alpha$-Calabi flow (\ref{combinatorial alpha Calabi flow}) exists for all time and converges exponentially fast to $\overline{u}$.
\end{description}
\end{theorem}
\begin{remark}\label{Calabi flow not extendable}
In the case of $\alpha=0$, Theorem \ref{main theorem 2} was proved by the first author in \cite{Xu 21}.
Similar to the case in \cite{Xu 21}, for generic initial discrete conformal structure in Definition \ref{edge length definition}, the combinatorial $\alpha$-Calabi flow can not be extended in the way used for the combinatorial $\alpha$-Ricci flow in this paper.
The global convergence of the combinatorial $\alpha$-Calabi flow (\ref{combinatorial alpha Calabi flow}) is not known up to now.
In the special case of vertex scaling, which corresponds to $\varepsilon \equiv 0$ for the discrete conformal structure in Definition \ref{edge length definition}, one can do surgery by flipping under Delaunay condition
on the combinatorial $\alpha$-Calabi flow to extend it and obtain the global convergence of the
combinatorial $\alpha$-Calabi flow with surgery. Please refer to \cite{Xu I,Zhu Xu,Xu Z} for more information on this.
\end{remark}

The paper is organized as follows. In Section \ref{section 2}, we study the rigidity for combinatorial $\alpha$-curvature of
Euclidean and hyperbolic discrete conformal structures
on polyhedral surfaces and prove a generalization of Theorem \ref{global rigidity theorem}.
In Section \ref{section 3}, we study the combinatorial Yamabe problem for combinatorial $\alpha$-curvature of discrete conformal structures using
combinatorial $\alpha$-Ricci flow and combinatorial $\alpha$-Calabi flow and prove generalizations of Theorem \ref{main theorem 1} and Theorem \ref{main theorem 2}.
In Section \ref{section 4}, we discuss some open problems related to combinatorial $\alpha$-curvatures and combinatorial $\alpha$-curvature flows.
\\
\\
\textbf{Acknowledgements}\\[8pt]
The paper is motivated by the ideas developed when the first author was visiting Professor Feng Luo at Rutgers University.
The first author thanks Professor Feng Luo for his invitation to Rutgers University and the communications during and after the visit.
Both authors thank the referees for their valuable suggestions to modify the paper substantially.
The research of the first author is supported by Fundamental Research Funds for the Central Universities under grant no. 2042020kf0199.

\section{Rigidity of combinatorial $\alpha$-curvatures}\label{section 2}

\subsection{Euclidean discrete conformal structures}
By Definition \ref{combinatorial alpha_curvature definition}, $u_i=f_i $ for all $i\in V$ in Euclidean background geometry.
Set $r_i=e^{u_i}$ for all $i\in V$, we also call $r\in \mathbb{R}^N_{>0}$ as a Euclidean discrete conformal structure.
The admissible space $\Omega^{E}_{ijk}$ of Euclidean discrete conformal structures
for a triangle $\{ijk\}\in F$ in $(M,\mathcal{T},\varepsilon,\eta)$  is defined to be
the set of $({r_i,r_j,r_k})\in \mathbb{R}^3_{>0}$ such that the triangle with edge lengths
$l_{ij}, l_{ik}, l_{jk}$ defined by (\ref{Euclidean edge length})  exists in 2-dimensional Euclidean space $\mathbb{E}^2$, i.e.
\begin{equation*}
\Omega^{E}_{ijk}=\{(r_i,r_j,r_k)\in \mathbb{R}_{>0}^3| l_{rs}+l_{rt}>l_{ts}, \{r,s,t\}=\{i,j,k\} \}.
\end{equation*}
The admissible space of Euclidean discrete conformal structures on $(M,\mathcal{T},\varepsilon,\eta)$ is
defined to be the vectors $r\in \mathbb{R}^N$ such that  $({r_i,r_j,r_k})\in \Omega^{E}_{ijk}$ for every triangle $\{ijk\}\in F$ and
we use $\Omega^{E}$ to denote it. One can also define the admissible space in terms of $f$, here we take the parameter $r$ for simplification of notations.

\begin{lemma}[\cite{Xu 21}]\label{Euclidean extended lemma}
Suppose $(M,\mathcal{T},\varepsilon,\eta)$ is a weighted triangulated connected closed surface with the weights $\varepsilon: V\rightarrow \{0,1\}$ and $\eta: E\rightarrow \mathbb{R}$ satisfying the structure conditions (\ref{structure condition 1}) and (\ref{structure condition 2}).
$\{ijk\}\in F$ is a topological triangle of the weighted triangulated surface $(M,\mathcal{T},\varepsilon,\eta)$.
\begin{description}
  \item[(1)] The admissible space $\Omega^{E}_{ijk}$ of Euclidean discrete conformal structures
  for $\{ijk\}\in F$  is nonempty and simply connected with analytic boundary.
  \item[(2)] The matrix $\Lambda^E_{ijk}:=\frac{\partial (\theta_{i}^{jk}, \theta_{j}^{ik}, \theta_{k}^{ij})}{\partial (u_i, u_j, u_k)}$ is symmetric and negative semi-definite with rank 2 and kernel $\{c(1,1,1)^T|c\in \mathbb{R} \}$ on $\Omega^{E}_{ijk}$, which implies that the matrix $\Lambda^E:=\frac{\partial(K_1,...,K_N)}{\partial(u_1,...,u_N)}$ is symmetric and positive semi-definite with rank $N-1$ and kernel $\{c(1,1,\cdots, 1)^T\in \mathbb{R}^N|c\in \mathbb{R} \}$ on $\Omega^E$.
  \item[(3)] The inner angles $\theta_i^{jk},\theta_j^{ik},\theta_k^{ij}$ defined for $(r_i,r_j,r_k)\in\Omega^{E}_{ijk}$ could be extended by constants to be continuous functions $\widetilde{\theta}_i^{jk},\widetilde{\theta}_j^{ik},\widetilde{\theta}_k^{ij}$ defined for  $(r_i,r_j,r_k)\in\mathbb{R}^3_{>0}$.
\end{description}
\end{lemma}

By Lemma \ref{Euclidean extended lemma}, we can extend the classical combinatorial curvature $K$ defined on $\Omega^{E}$ to be defined for $r\in \mathbb{R}^N_{>0}$ by setting
\begin{equation*}
\widetilde{K}_i=2\pi-\sum_{\{ijk\}\in F}\widetilde{\theta}_i^{jk},
\end{equation*}
which still satisfies the discrete Gauss-Bonnet formula $\sum_{i=1}^N\widetilde{K}_i=2\pi\chi(M)$.
As a result, the combinatorial $\alpha$-curvature $R_{\alpha}$ could be extended by setting $$\widetilde{R}_{\alpha,i}=\frac{\widetilde{K}_i}{e^{\alpha u_i}}=\frac{\widetilde{K}_i}{r_i^\alpha}$$
for any $r\in \mathbb{R}^N_{>0}$ and $i\in V$.
We call both vectors $u=(u_1,...,u_N)\in \mathbb{R}^N$ and $r=(r_1,...,r_N)=(e^{u_1},...,e^{u_N})\in \mathbb{R}^N_{>0}$ as  generalized Euclidean discrete conformal structures on $(M,\mathcal{T},\varepsilon,\eta)$, if there is no confusion.

According to Lemma \ref{Euclidean extended lemma}, the Ricci energy function
\begin{equation*}
F_{ijk}(u_i,u_j,u_k)=\int^{(u_i,u_j,u_k)}_{0}\theta_i^{jk}du_i+\theta_j^{ik}du_j+\theta_k^{ij}du_k
\end{equation*}
for a triangle $\{ijk\}\in F$ is well-defined on $\mathcal{U}^{E}_{ijk}=\ln\Omega^{E}_{ijk}$. Furthermore, $F_{ijk}(u_i,u_j,u_k)$ is locally concave on $\mathcal{U}^{E}_{ijk}$ and locally strictly concave on $\mathcal{U}^{E}_{ijk}\cap \{u_i+u_j+u_k=0\}$ with $\nabla_{u_i}F_{ijk}=\theta_i^{jk}$ and $F_{ijk}(u_i+t,u_j+t,u_k+t)=F_{ijk}(u_i,u_j,u_k)+t\pi$.

For further applications, we need to extend $F_{ijk}$ to be a globally defined function.
Let us recall Luo's generalization \cite{Luo3} of Bobenko-Pinkall-Spingborn's extension \cite{BPS} first.
\begin{definition}[\cite{Luo3}, Definition 2.3]
A differential 1-form $w=\sum_{i=1}^n a_i(x)dx^i$ in an open subset $U\subset \mathbb{R}^n$ is said to be continuous if each $a_i(x)$ is continuous on $U$. A continuous differential 1-form $w$ is called closed if $\int_{\partial \tau}w=0$ for each triangle $\tau\subset U$.
\end{definition}
\begin{theorem}[\cite{Luo3}, Corollary 2.6]\label{Luo's convex extention}
Suppose $X\subset \mathbb{R}^n$ is an open convex set and $A\subset X$ is an open subset of $X$ bounded by a real analytic codimension-1 submanifold in $X$. If $w=\sum_{i=1}^na_i(x)dx_i$ is a continuous closed 1-form on $A$ so that $F(x)=\int_a^x w$ is locally convex on $A$ and each $a_i$ can be extended continuous to $X$ by constant functions to a function $\widetilde{a}_i$ on $X$, then $\widetilde{F}(x)=\int_a^x\sum_{i=1}^n\widetilde{a}_i(x)dx_i$ is a $C^1$-smooth
convex function on $X$ extending $F$.
\end{theorem}

Combining Lemma \ref{Euclidean extended lemma} and Theorem \ref{Luo's convex extention}, $F_{ijk}(u_i,u_j,u_k)$ defined on $\mathcal{U}^{E}_{ijk}$ could be extended to be
\begin{equation}\label{Euclidean extension of Fijk}
\widetilde{F}_{ijk}(u_i,u_j,u_k)=\int^{(u_i,u_j,u_k)}_{0}\widetilde{\theta}_i^{jk}du_i+\widetilde{\theta}_j^{ik}du_j
+\widetilde{\theta}_k^{ij}du_k,
\end{equation}
which is a $C^1$-smooth concave function defined on $\mathbb{R}^3$ with $\nabla_{u_i}\widetilde{F}_{ijk}=\widetilde{\theta}_i^{jk}$.
Using the extension $\widetilde{F}_{ijk}$ of Ricci energy function $F_{ijk}$,
we can prove the following rigidity for the extended combinatorial $\alpha$-curvature $\widetilde{R}_{\alpha}$,
which is a generalization of Theorem \ref{global rigidity theorem} (a).

\begin{theorem}\label{Euclidean global rigidity theorem}
Suppose $(M,\mathcal{T},\varepsilon,\eta)$ is a weighted triangulated  connected closed surface with the weights $\varepsilon: V\rightarrow \{0,1\}$ and $\eta: E\rightarrow \mathbb{R}$ satisfying the structure conditions (\ref{structure condition 1}) and (\ref{structure condition 2}). $\alpha\in \mathbb{R}$ is a constant and $\overline{R}: V\rightarrow \mathbb{R}$ is a given function.
Suppose there exist $r_A\in \Omega^{E}$ and $r_B\in \mathbb{R}^N_{>0}$ with the same extended combinatorial $\alpha$-curvature $\overline{R}$.
If $\alpha\overline{R}\equiv 0$, then $r_A=\lambda r_B$ for some positive constant $\lambda\in \mathbb{R}_{>0}$. If $\alpha\overline{R}\leq 0$ and $\alpha\overline{R}\not\equiv 0$, then $r_A=r_B$.
\end{theorem}
\proof
Define the following Ricci energy $F(u)$ for $\overline{R}$ by
\begin{equation*}
F(u)=-\sum_{\{ijk\}\in F}F_{ijk}(u_i,u_j,u_k)
+\int^u_{u_0}\sum_{i=1}^{N}(2\pi-\overline{R}_ir^{\alpha}_i)du_i.
\end{equation*}
By direct calculations, we have
\begin{equation*}
\nabla_{u_i}F(u)=-\sum_{\{ijk\}\in F}\theta_i^{jk}+2\pi-\overline{R}_ir^{\alpha}_i
=K_i-\overline{R}_ir^{\alpha}_i
\end{equation*}
and
\begin{equation*}
\mathrm{Hess}_u F=\Lambda^E-\alpha \left(
 \begin{array}{ccc}
  \overline{R}_1r^{\alpha}_1  &  &\\
     & \ddots &  \\
      & & \overline{R}_Nr^{\alpha}_N \\
  \end{array}
 \right)
\end{equation*}
for $u\in \mathcal{U}^E=\ln \Omega^{E}$.
By Lemma \ref{Euclidean extended lemma}, if $\alpha\overline{R}\equiv 0$, then $\mathrm{Hess}_u F$ is positive semi-definite with kernel $\{c(1,1,\cdots, 1)^T|c\in \mathbb{R} \}$ and $F$ is locally convex on $\mathcal{U}^E$.
If $\alpha\overline{R}\leq 0$ and $\alpha\overline{R}\not\equiv 0$, then $\mathrm{Hess}_u F$ is positive definite and $F$ is locally strictly convex on $\mathcal{U}^E$.

By the extension $\widetilde{F}_{ijk}(u_i,u_j,u_k)$ of $F_{ijk}(u_i,u_j,u_k)$ in (\ref{Euclidean extension of Fijk}),
the Ricci energy function $F(u)$ defined on $\mathcal{U}^E$ could be extended to be
\begin{equation}\label{Euclidean extended F}
\widetilde{F}(u)=-\sum_{\{ijk\}\in F}\widetilde{F}_{ijk}(u_i,u_j,u_k)
+\int^u_{u_0}\sum_{i=1}^{N}(2\pi-\overline{R}_ir^{\alpha}_i)du_i,
\end{equation}
which is a $C^1$-smooth convex function defined on $\mathbb{R}^N$. Furthermore,
\begin{equation*}
\nabla_{u_i}\widetilde{F}=-\sum_{\{ijk\}\in F}\widetilde{\theta}_i^{jk}+2\pi-\overline{R}_ir^{\alpha}_i=\widetilde{K}_i-\overline{R}_ir^{\alpha}_i,
\end{equation*}
where $\widetilde{K}_i=2\pi-\sum_{\{ijk\}\in F}\widetilde{\theta}_i^{jk}$.
Then $\widetilde{F}(u)$ is convex on $\mathbb{R}^N$ and
locally strictly convex on $\mathcal{U}^E\cap \{\sum^N_{i=1}u_i=0\}$
in the case of $\alpha\overline{R}\equiv 0$. Similarly, $\widetilde{F}(u)$ is convex on $\mathbb{R}^N$ and locally strictly convex on $\mathcal{U}^E$ in the case of $\alpha\overline{R}\leq 0$ and $\alpha\overline{R}\not\equiv 0$.

If there exist $r_A\in \Omega^{E}$ and $r_B\in \mathbb{R}^N_{>0}$ with the same extended combinatorial $\alpha$-curvature $\overline{R}$, set $h(t)=\widetilde{F}((1-t)u_A+tu_B),\ t\in [0,1].$
Then $h(t)$ is a $C^1$ convex function with
$$h'(t)
=\sum_{i=1}^N\nabla_{u_i} \widetilde{F}|_{(1-t)u_A+tu_B}\cdot (u_{B,i}-u_{A,i})
=\sum_{i=1}^N(\widetilde{K}_i-\overline{R}_ir^{\alpha}_i)|_{(1-t)u_A+tu_B}\cdot (u_{B,i}-u_{A,i}).$$
By the assumption that $\widetilde{R}_\alpha(u_A)=\widetilde{R}_\alpha(u_B)=\overline{R}$, we have $h'(0)=h'(1)=0$,
which implies $h'(t)\equiv0$ by the convexity of $h(t)$.
Note that $u_A$ is in the open subset $\mathcal{U}^E\subseteq\mathbb{R}^N$, there exists $0<\epsilon<1$ such that $(1-t)u_A+tu_B\in \mathcal{U}^E$ and $h(t)$ is smooth for $t\in[0,\epsilon)$.

In the case of $\alpha\overline{R}\leq 0$ and $\alpha\overline{R}\not\equiv 0$, $\widetilde{F}(u)$ is locally strictly convex on $\mathcal{U}^E$, which implies $h(t)$ is locally strictly convex on $[0,\epsilon)$ and $h'(\epsilon)>0$ by $h'(0)=0$. This contradicts $h'(t)\equiv0$ on $[0,1]$. So $u_A=u_B$, which implies $r_A=r_B$.

In the case of $\alpha\overline{R}\equiv 0$, we have $h(t)$ is $C^1$ convex on $[0,1]$ and smooth on $[0,\epsilon)$.  $h'(t)\equiv0$ for $t\in [0, 1]$ implies that $h''(t)\equiv0$ for $t\in [0, \epsilon)$.
Note that
$$h''(t)=(u_B-u_A)^T\cdot \mathrm{Hess} F|_{(1-t)u_A+tu_B} \cdot (u_B-u_A)$$
for $t\in [0, \epsilon)$ and $\mathrm{Hess} F$ is positive semi-definite with kernel $\{c(1,1,\cdots, 1)^T\in \mathbb{R}^N|c\in \mathbb{R} \}$,
we have $u_B-u_A=\mu(1,1,\cdots, 1)^T$ for some constant $\mu\in \mathbb{R}$, which implies $r_B=\lambda r_A$ with $\lambda=e^\mu$.
\qed

As a corollary of Theorem \ref{Euclidean global rigidity theorem}, we have the following result on the global rigidity for discrete conformal structures with constant combinatorial $\alpha$-curvature.
\begin{corollary}\label{Euclidean global rigidity corollary}
Suppose $(M,\mathcal{T},\varepsilon,\eta)$ is a weighted triangulated  connected closed surface with the weights $\varepsilon: V\rightarrow \{0,1\}$ and $\eta: E\rightarrow \mathbb{R}$ satisfying the structure conditions (\ref{structure condition 1}) and (\ref{structure condition 2}). $\alpha\in \mathbb{R}$ is a constant. If $\alpha\chi(M)= 0$, there exists at most one Euclidean discrete conformal structure $f: V\rightarrow \mathbb{R}$ with constant combinatorial $\alpha$-curvature up to a vector $c(1,1,...,1)^T,\ c\in \mathbb{R}$. If $\alpha\chi(M)<0$, there exists at most one Euclidean discrete conformal structure $f: V\rightarrow \mathbb{R}$ with constant combinatorial $\alpha$-curvature.
\end{corollary}

\subsection{Hyperbolic discrete conformal structures}
The admissible space $\Omega^{H}_{ijk}$ of hyperbolic discrete conformal structures for a triangle $\{ijk\}\in F$
is defined to be the set of $({f_i,f_j,f_k})\in \mathbb{R}^3$ such that the triangle with edge lengths
$l_{ij}, l_{ik}, l_{jk}$ defined by (\ref{hyperbolic edge length})
 exists in 2-dimensional hyperbolic space $\mathbb{H}^2$, i.e.
\begin{equation*}
\Omega^{H}_{ijk}=\{(f_i,f_j,f_k)\in \mathbb{R}^3| l_{rs}+l_{rt}>l_{ts}, \{r,s,t\}=\{i,j,k\} \}.
\end{equation*}
The admissible space of hyperbolic discrete conformal structures on $(M,\mathcal{T},\varepsilon,\eta)$ is
defined to be the vectors $f\in \mathbb{R}^N$ such that  $({f_i,f_j,f_k})\in \Omega^{H}_{ijk}$ for every triangle $\{ijk\}\in F$ and
we use $\Omega^{H}$ to denote it.

In the hyperbolic background geometry,
for the map $u$ defined by (\ref{hyperblic DCS in u}), we have
 $u=(u_1,..,u_N)\in \mathbb{R}^{N_0}\times \mathbb{R}^{N_1}_{<0}$ by Remark \ref{remark on u}, where $N_0$ is the number of the vertices $i\in V$ with $\varepsilon_i=0$ and $N_1=N-N_0$. For simplicity, we also call $u\in \mathbb{R}^{N_0}\times \mathbb{R}^{N_1}_{<0}$ as hyperbolic discrete conformal structure.

\begin{lemma}[\cite{Xu 21}]\label{hyperbolic extended lemma}
Suppose $(M,\mathcal{T},\varepsilon,\eta)$ is a weighted triangulated  connected closed surface with the weights $\varepsilon: V\rightarrow \{0,1\}$ and $\eta: E\rightarrow \mathbb{R}$ satisfying the structure conditions (\ref{structure condition 1}) and (\ref{structure condition 2}).
$\{ijk\}\in F$ is a topological triangle of the weighted triangulated surface $(M,\mathcal{T},\varepsilon,\eta)$.
\begin{description}
  \item[(1)] The admissible space $\Omega^{H}_{ijk}$ of hyperbolic discrete conformal structures
  for $\{ijk\}\in F$  is nonempty and simply connected with analytic boundary.
  \item[(2)] The matrix $\Lambda^H_{ijk}:=\frac{\partial (\theta_{i}^{jk}, \theta_{j}^{ik}, \theta_{k}^{ij})}{\partial (u_i, u_j, u_k)}$ is symmetric and negative definite on $\Omega^{H}_{ijk}$, which implies the matrix $\Lambda^H:=\frac{\partial(K_1,...,K_N)}{\partial(u_1,...,u_N)}$ is symmetric and positive definite on $\Omega^{H}$.
  \item[(3)] The inner angles $\theta_i^{jk},\theta_j^{ik},\theta_k^{ij}$ defined for $(f_i,f_j,f_k)\in\Omega^H_{ijk}$ could be extended by constants to be continuous functions $\widetilde{\theta}_i^{jk},\widetilde{\theta}_j^{ik},\widetilde{\theta}_k^{ij}$ defined for $(f_i,f_j,f_k)\in\mathbb{R}^3$.
\end{description}
\end{lemma}

By Lemma \ref{hyperbolic extended lemma}, we can extend the classical combinatorial curvature $K$ defined on $\Omega^{H}$ to be defined on $\mathbb{R}^N$ by setting
$\widetilde{K}_i=2\pi-\sum_{\{ijk\}\in F}\widetilde{\theta}_i^{jk}$.
As a result, the combinatorial $\alpha$-curvature $R_\alpha$ defined for $f\in \Omega^{H}$ could be extended to be defined for  $f\in \mathbb{R}^N$ by setting $\widetilde{R}_{\alpha,i}=\frac{\widetilde{K}_i}{e^{\alpha u_i}}$.
We call both vectors $f=(f_1,...,f_N)\in \mathbb{R}^N$ and $u=(u_1,...,u_N)\in \mathbb{R}^{N_0}\times \mathbb{R}^{N_1}_{<0}$ as generalized hyperbolic discrete conformal structures on $(M,\mathcal{T},\varepsilon,\eta)$, if there is no confusion.

By Lemma \ref{hyperbolic extended lemma}, the Ricci energy function
\begin{equation*}
F_{ijk}(u_i,u_j,u_k)=\int^{(u_i,u_j,u_k)}_{0}\theta_i^{jk}du_i+\theta_j^{ik}du_j+\theta_k^{ij}du_k
\end{equation*}
for a triangle $\{ijk\}\in F$ is a well-defined locally strictly concave function defined on $\mathcal{U}^H_{ijk}$ with $\nabla_{u_i}F_{ijk}=\theta_i^{jk}$, where $\mathcal{U}^H_{ijk}$ is the image of $\Omega^{H}_{ijk}$ under the map (\ref{hyperblic DCS in u}).
Combining Lemma \ref{hyperbolic extended lemma} and Theorem \ref{Luo's convex extention}, $F_{ijk}(u_i,u_j,u_k)$ defined on $\mathcal{U}^{H}_{ijk}$ could be extended to be
\begin{equation*}
\widetilde{F}_{ijk}(u_i,u_j,u_k)=\int^{(u_i,u_j,u_k)}_{0}\widetilde{\theta}_i^{jk}du_i+\widetilde{\theta}_j^{ik}du_j
+\widetilde{\theta}_k^{ij}du_k,
\end{equation*}
which is a $C^1$-smooth concave function defined on $\mathbb{R}^{n_0}\times \mathbb{R}^{n_1}_{<0}$ with $\nabla_{u_i}\widetilde{F}_{ijk}=\widetilde{\theta}_i$, where $n_0$ is the number of the vertices in $\{i,j,k\}$ with $\varepsilon=0$ and $n_1=3-n_0$.

Paralleling to Theorem \ref{Euclidean global rigidity theorem} for generalized Euclidean discrete conformal structures, we prove the following rigidity of the extended combinatorial $\alpha$-curvature $\widetilde{R}_\alpha$ for generalized hyperbolic discrete conformal structures, which is a generalization of Theorem \ref{global rigidity theorem} (b).

\begin{theorem}\label{hyperbolic global rigidity theorem}
Suppose $(M,\mathcal{T},\varepsilon,\eta)$ is a weighted triangulated connected closed surface with the weights $\varepsilon: V\rightarrow \{0,1\}$ and $\eta: E\rightarrow \mathbb{R}$ satisfying the structure conditions (\ref{structure condition 1}) and (\ref{structure condition 2}). $\alpha\in \mathbb{R}$ is a constant and $\overline{R}: V\rightarrow \mathbb{R}$ is a given function with $\alpha\overline{R}\leq 0$.
If there exist $f_A\in \Omega^{H}$ and $f_B\in \mathbb{R}^N$ with the same extended combinatorial $\alpha$-curvature $\overline{R}$, then $f_A=f_B$.
\end{theorem}
\proof
Define the following Ricci energy $F(u)$ for $\overline{R}$ by
\begin{equation*}
F(u)=-\sum_{\{ijk\}\in F}F_{ijk}(u_i,u_j,u_k)
+\int^u_{u_0}\sum_{i=1}^{N}(2\pi-\overline{R}_ie^{\alpha u_i})du_i.
\end{equation*}
Then
\begin{equation*}
\mathrm{Hess}_u F(u)=\Lambda^H-\alpha \left(
 \begin{array}{ccc}
  \overline{R}_1e^{\alpha u_1}  &  &\\
     & \ddots &  \\
      & & \overline{R}_Ne^{\alpha u_N} \\
  \end{array}
 \right)
\end{equation*}
on $\mathcal{U}^H=u(\Omega^H)$.
Combining Lemma \ref{hyperbolic extended lemma} with the condition $\alpha\overline{R}\leq 0$, $\mathrm{Hess}_u F$ is positive definite and $F$ is locally strictly convex on $\mathcal{U}^H$.

Parallelling to the Euclidean case,
the Ricci energy function $F(u)$ defined on $\mathcal{U}^H$ could be extended to be
\begin{equation}\label{hyperbolic extended F}
\widetilde{F}(u)=-\sum_{\{ijk\}\in F}\widetilde{F}_{ijk}(u_i,u_j,u_k)
+\int^u_{u_0}\sum_{i=1}^{N}(2\pi-\overline{R}_ie^{\alpha u_i})du_i,
\end{equation}
which is a $C^1$-smooth convex function defined on $\mathbb{R}^{N_0}\times \mathbb{R}^{N_1}_{<0}$ and locally strictly convex on $\mathcal{U}^{H}\subset \mathbb{R}^{N_0}\times \mathbb{R}^{N_1}_{<0}$ with
\begin{equation*}
\nabla_{u_i}\widetilde{F}
=-\sum_{\{ijk\}\in F}\widetilde{\theta}_i+2\pi-\overline{R}_ie^{\alpha u_i}
=\widetilde{K}_i-\overline{R}_ie^{\alpha u_i}.
\end{equation*}

If there exist $f_A\in \Omega^{H}$ and $f_B\in \mathbb{R}^N$ with the same extended combinatorial $\alpha$-curvature $\overline{R}$, set $h(t)=\widetilde{F}((1-t)u_A+tu_B),\ t\in [0,1]$,
where $u_A=u(f_A), u_B=u(f_B)$ with $u$ given by (\ref{hyperblic DCS in u}).
Then $h(t)$ is a $C^1$ convex function for $t\in [0,1]$.
The proof in the following is parallelling to that for Theorem \ref{Euclidean global rigidity theorem}, we omit the details here.
\qed

\section{Deformation of discrete conformal structures}\label{section 3}
In this section, we use combinatorial curvature flows to study the combinatorial Yamabe problem for
combinatorial $\alpha$-curvature
of discrete conformal structures on polyhedral surfaces.
In the Euclidean background geometry, we study constant combinatorial $\alpha$-curvature problem
and prescribed combinatorial $\alpha$-curvature problem using the combinatorial $\alpha$-Ricci flow (\ref{normalized combinatorial alpha Ricci flow}), (\ref{modified combinatorial alpha Ricci flow}) and the combinatorial $\alpha$-Calabi flow (\ref{combinatorial alpha Calabi flow}).
In hyperbolic background geometry,
we take the constant combinatorial $\alpha$-curvature problem as a special case of
prescribed combinatorial $\alpha$-curvature problem
and use the modified combinatorial $\alpha$-Ricci flow (\ref{modified combinatorial alpha Ricci flow})
and the following modified combinatorial $\alpha$-Calabi flow
\begin{equation}\label{modified combinatorial alpha Calabi flow}
\frac{du_i}{dt}=\Delta_\alpha(R_{\alpha}-\overline{R})_i
\end{equation}
to study the prescribed combinatorial $\alpha$-curvature problem.

\subsection{Local convergence of combinatorial $\alpha$-curvature flows}

We have the following properties of the normalized combinatorial $\alpha$-Ricci flow (\ref{normalized combinatorial alpha Ricci flow}) and the combinatorial $\alpha$-Calabi flow (\ref{combinatorial alpha Calabi flow}) in Euclidean background geometry.

\begin{lemma}\label{invariant along flow}
In Euclidean background geometry, if $\alpha=0$, then $\sum^N_{i=1}u_i$ is invariant along the normalized combinatorial $\alpha$-Ricci flow (\ref{normalized combinatorial alpha Ricci flow}) and the combinatorial $\alpha$-Calabi flow (\ref{combinatorial alpha Calabi flow}). If $\alpha\neq0$, $||r||^\alpha_\alpha=\sum^N_{i=1}r^\alpha_i=\sum^N_{i=1}e^{\alpha u_i}$ is invariant along the normalized combinatorial $\alpha$-Ricci flow (\ref{normalized combinatorial alpha Ricci flow}) and the combinatorial $\alpha$-Calabi flow (\ref{combinatorial alpha Calabi flow}).
\end{lemma}
\proof
The case of $\alpha=0$ has been proved in \cite{Xu 21}, we only prove the case of $\alpha\neq 0$ here.
If $\alpha\neq0$, along the equation (\ref{normalized combinatorial alpha Ricci flow}), we have
\begin{equation*}
\frac{d(\sum^N_{i=1}e^{\alpha u_i})}{dt}
=\sum^N_{i=1}\alpha e^{\alpha u_i}\frac{du_i}{dt}
=\alpha\left(2\pi\chi(M)-\sum^N_{i=1}K_i\right)
=0,
\end{equation*}
which implies $\sum^N_{i=1}e^{\alpha u_i}$ is invariant along the normalized combinatorial $\alpha$-Ricci flow (\ref{normalized combinatorial alpha Ricci flow}).
Similarly, along the equation (\ref{combinatorial alpha Calabi flow}), we have
\begin{equation*}
\frac{d(\sum^N_{i=1}e^{\alpha u_i})}{dt}
=-\alpha\sum^N_{j=1}\sum^N_{i=1}(\Lambda^E)_{ij}R_{\alpha,j}
=0
\end{equation*}
by Lemma \ref{Euclidean extended lemma},
which implies $\sum^N_{i=1}e^{\alpha u_i}$ is invariant along the combinatorial $\alpha$-Calabi flow (\ref{combinatorial alpha Calabi flow}).
\qed

The normalized combinatorial $\alpha$-Ricci flow (\ref{normalized combinatorial alpha Ricci flow})
and the combinatorial $\alpha$-Calabi flow (\ref{combinatorial alpha Calabi flow}) are ODE systems
with smooth coefficients.
Therefore, the solutions always exist locally around the initial time $t=0$.
We further have the following result on the longtime behaviors for the solutions of normalized combinatorial $\alpha$-Ricci flow (\ref{normalized combinatorial alpha Ricci flow}) and combinatorial $\alpha$-Calabi flow (\ref{combinatorial alpha Calabi flow}) in  Euclidean background geometry, which is a generalization of Theorem \ref{main theorem 2} (a).

\begin{theorem}\label{alpha flows with small initial energy 1}
Suppose $(M,\mathcal{T},\varepsilon,\eta)$ is a weighted triangulated connected closed surface with the weights $\varepsilon: V\rightarrow \{0,1\}$ and $\eta: E\rightarrow \mathbb{R}$ satisfying the structure conditions (\ref{structure condition 1}) and (\ref{structure condition 2}). $\alpha\in \mathbb{R}$ is a constant.
\begin{description}
  \item[(1)] If the solution of normalized combinatorial $\alpha$-Ricci flow (\ref{normalized combinatorial alpha Ricci flow}) or the solution of combinatorial $\alpha$-Calabi flow (\ref{combinatorial alpha Calabi flow}) converges in $\mathcal{U}^E$, there exists a discrete conformal structure with constant combinatorial $\alpha$-curvature.
  \item[(2)] Suppose that there exists a discrete conformal structure $\overline{u}\in \mathcal{U}^E$ with constant combinatorial $\alpha$-curvature and $\alpha \chi(M)\leq 0$, then there exists a constant $\delta >0$ such that if
      the initial value $u(0)$ satisfies
      $||R_\alpha(u(0))-R_\alpha(\overline{u})||<\delta$ and
$\sum^N_{i=1}e^{\alpha u(0)}=\sum^N_{i=1}e^{\alpha\overline{u}_i}$ in the case of $\alpha\neq 0$ or $||u(0)-\overline{u}||<\delta$ and
 $\sum^N_{i=1}u_i(0)=\sum^N_{i=1}\overline{u}_i$ in the case of $\alpha=0$,
      the solutions of normalized combinatorial $\alpha$-Ricci flow (\ref{normalized combinatorial alpha Ricci flow}) and the combinatorial $\alpha$-Calabi flow (\ref{combinatorial alpha Calabi flow}) exist for all time and converge exponentially fast to $\overline{u}$ respectively.
\end{description}
\end{theorem}
\proof
The case of $\alpha=0$ has been proved in \cite{Xu 21}, we only prove the case of $\alpha\neq 0$ here.
Suppose $u(t)$ is a solution of the normalized combinatorial $\alpha$-Ricci flow (\ref{normalized combinatorial alpha Ricci flow}). If $\overline{u}:=u(+\infty)=\lim_{t\rightarrow +\infty} u(t)$ exists in $\mathcal{U}^E$, then $R_\alpha(\overline{u})=\lim_{t\rightarrow +\infty}R_\alpha(u(t))$ exists.
Furthermore, there exists a sequence $\xi_n\in(n,n+1)$ such that
$$u_i(n+1)-u_i(n)=u'_i(\xi_n)=R_{\alpha,av}-R_{\alpha,i}(u(\xi_n))\rightarrow 0\  \text{as} \ n\rightarrow +\infty,$$
which implies $R_{\alpha}(\overline{u})=R_{\alpha,av}$ and $\overline{u}$ is a discrete conformal structure with constant combinatorial $\alpha$-curvature $R_{\alpha,av}$.
Similarly, if the combinatorial $\alpha$-Calabi flow (\ref{combinatorial alpha Calabi flow}) converges, there exists a discrete conformal structure with constant combinatorial $\alpha$-curvature.

Suppose there exists a discrete conformal structure $\overline{u}\in \mathcal{U}^E$ with constant combinatorial $\alpha$-curvature $R_{\alpha,av}$. For the normalized combinatorial $\alpha$-Ricci flow (\ref{normalized combinatorial alpha Ricci flow}), set $\Gamma_i(u)=R_{\alpha,av}-R_{\alpha,i}$. By direct calculations, we have
\begin{equation*}
\begin{aligned}
D\Gamma|_{u=\overline{u}}
&=\alpha R_{\alpha,av}I-\Sigma^{-\alpha}\left(\Lambda^E+\alpha R_{\alpha,av}\frac{r^\alpha\cdot (r^\alpha)^T}{||r||^\alpha_\alpha}\right)\\
&=-\Sigma^{-\frac{\alpha}{2}}\left(L-\alpha R_{\alpha,av}\left[I-\frac{r^\alpha\cdot (r^\alpha)^T}{||r||^\alpha_\alpha}\right]\right)\Sigma^{\frac{\alpha}{2}},
\end{aligned}
\end{equation*}
where $\Sigma=\text{diag}\{r_1,r_2,...,r_N\}$, $L=\Sigma^{-\frac{\alpha}{2}}\Lambda^E\Sigma^{-\frac{\alpha}{2}}$.
Note that the matrix $I-\frac{r^\alpha\cdot (r^\alpha)^T}{||r||^\alpha_\alpha}$ has eigenvalues 1 ($N-1$ times) and 0 (1 time) and its kernel is $\{cr^{\frac{\alpha}{2}}|c\in \mathbb{R}\}$.
Further note that the matrix $L$ is positive semi-definite with 1-dimensional kernel $\{cr^{\frac{\alpha}{2}}|c\in \mathbb{R}\}$ by Lemma \ref{Euclidean extended lemma}.
Therefore, by Lemma \ref{invariant along flow} and the condition $\alpha \chi(M)\leq 0$, $D\Gamma|_{u=\overline{u}}$ has $N-1$ negative eigenvalues and a zero eigenvalue with 1-dimensional kernel $\{c(1,1,...,1)\}$.
Let $\widetilde{u}_i=e^{\frac{\alpha}{2}\overline{u}_i}u_i
=\overline{r}^{\frac{\alpha}{2}}_iu_i$ and $\overline{\Sigma}=\Sigma(\overline{u})=\text{diag}
\{\overline{r}_1,\overline{r}_2,...,\overline{r}_N\}$,
then $\widetilde{u}=\overline{\Sigma}^{\frac{\alpha}{2}}u$.
The combinatorial $\alpha$-Ricci flow with respect to $\widetilde{u}$ is
$\frac{d\widetilde{u}_i}{dt}
=\overline{r}_i^{\frac{\alpha}{2}}\frac{du_i}{dt}
=\overline{r}_i^{\frac{\alpha}{2}}\Gamma_i(u)$.
Let $\widetilde{\Gamma}(\widetilde{u})
=\overline{\Sigma}^{\frac{\alpha}{2}}
\Gamma(\overline{\Sigma}^{-\frac{\alpha}{2}}\widetilde{u})$,
then the matrix
$$\left(\frac{\partial \widetilde{\Gamma}}{\partial \widetilde{u}}\right)\bigg|_{\widetilde{u}=
\overline{\Sigma}^{\frac{\alpha}{2}}\overline{u}}
=\overline{\Sigma}^{\frac{\alpha}{2}}\cdot\left(\frac{\partial \Gamma}{\partial u}\right)\bigg|_{u=\overline{u}}\cdot\overline{\Sigma}^{-\frac{\alpha}{2}}.$$
is symmetric and negative semi-definite with 1-dimensional kernel $\{c\overline{r}^{\frac{\alpha}{2}}|c\in \mathbb{R}\}$ by the property of $D\Gamma|_{u=\overline{u}}$.
The hypersurface $\{r\in \mathbb{R}^N|\sum^N_{i=1}r^\alpha_i=\sum^N_{i=1}e^{\alpha u_i}=\sum^N_{i=1}r^\alpha_i(0)\}$ is equivalent to $\{\widetilde{u}_i\in \mathbb{R}|\sum^N_{i=1}e^{\alpha e^{-\frac{\alpha}{2}\overline{u}_i}\widetilde{u}_i}=\sum^N_{i=1}r^\alpha_i(0)\}$,
the normal of which at $\widetilde{u}_i$ is
$n_i=\alpha e^{-\frac{\alpha}{2}\overline{u}_i}e^{\alpha e^{-\frac{\alpha}{2}\overline{u}_i}\widetilde{u}_i}=\alpha e^{-\frac{\alpha}{2}\overline{u}_i}e^{\alpha\overline{u}_i}
=\alpha e^{\frac{\alpha}{2}\overline{u}_i}$.
Thus the matrix $\left(\frac{\partial \widetilde{\Gamma}}{\partial \widetilde{u}}\right)\big|_{\widetilde{u}=
\overline{\Sigma}^{\frac{\alpha}{2}}\overline{u}}$ restricted to the hypersurface $\{r\in \mathbb{R}^N|\sum^N_{i=1}r^\alpha_i=\sum^N_{i=1}r^\alpha_i(0)\}$ is negative definite, which implies that $D\Gamma|_{u=\overline{u}}$ restricted to the hypersurface $\{r\in \mathbb{R}^N|\sum^N_{i=1}r^\alpha_i=\sum^N_{i=1}r^\alpha_i(0)\}$ is negative definite.
Therefore, $\overline{u}$ is a local attractor of the normalized combinatorial $\alpha$-Ricci flow (\ref{normalized combinatorial alpha Ricci flow}). Then the conclusion follows from Lyapunov Stability Theorem (\cite{Pontryagin}, Chapter 5).

Similarly, for the combinatorial $\alpha$-Calabi flow (\ref{modified combinatorial alpha Calabi flow}), set $\Gamma_i(u)=\Delta_\alpha R_{\alpha,i}$. By direct calculations, we have
\begin{equation*}
\begin{aligned}
D\Gamma|_{u=\overline{u}}
&=-\Sigma^{-\alpha}\Lambda^E\Sigma^{-\alpha}\Lambda^E+\alpha R_{\alpha,av}\Sigma^{-\alpha}\Lambda^E\\
&=-\Sigma^{-\frac{\alpha}{2}}\left(\Sigma^{-\frac{\alpha}{2}}\Lambda^E\Sigma^{-\alpha}\Lambda^E\Sigma^{-\frac{\alpha}{2}}
-\alpha R_{\alpha,av}\Sigma^{-\frac{\alpha}{2}}\Lambda^E\Sigma^{-\frac{\alpha}{2}}\right)\Sigma^{\frac{\alpha}{2}}.
\end{aligned}
\end{equation*}
By Lemma \ref{Euclidean extended lemma} and the condition $\alpha \chi(M)\leq 0$, $D\Gamma|_{u=\overline{u}}$ has $N-1$ negative eigenvalues and a zero eigenvalue with 1-dimensional kernel $\{c(1,1,...,1)\}$.
Using the above trick again, $D\Gamma|_{u=\overline{u}}$ restricted to the hypersurface $\{r\in \mathbb{R}^N|\sum^N_{i=1}r^\alpha_i=\sum^N_{i=1}r^\alpha_i(0)\}$ is negative definite. By Lemma \ref{invariant along flow},  $\overline{u}$ is a local attractor of the combinatorial $\alpha$-Calabi flow (\ref{combinatorial alpha Calabi flow}). Then the conclusion follows from Lyapunov Stability Theorem (\cite{Pontryagin}, Chapter 5).
\qed

Similar to Theorem \ref{alpha flows with small initial energy 1} for Euclidean background geometry, we have the following result on the longtime existence and convergence for the solutions of modified combinatorial $\alpha$-Ricci flow (\ref{modified combinatorial alpha Ricci flow}) and modified combinatorial $\alpha$-Calabi flow (\ref{modified combinatorial alpha Calabi flow}) in the hyperbolic background geometry, which is a generalization of Theorem \ref{main theorem 2} (b).

\begin{theorem}\label{alpha flows with small initial energy 2}
Suppose $(M,\mathcal{T},\varepsilon,\eta)$ is a weighted triangulated connected closed surface with the weights $\varepsilon: V\rightarrow \{0,1\}$ and $\eta: E\rightarrow \mathbb{R}$ satisfying the structure conditions (\ref{structure condition 1}) and (\ref{structure condition 2}). $\alpha\in \mathbb{R}$ is a constant and $\overline{R}: V\rightarrow \mathbb{R}$ is a given function.
\begin{description}
  \item[(1)] If the solution of the modified combinatorial $\alpha$-Ricci flow (\ref{modified combinatorial alpha Ricci flow}) or the solution of modified combinatorial $\alpha$-Calabi flow (\ref{modified combinatorial alpha Calabi flow}) converges in $\mathcal{U}^H$, there exists a hyperbolic discrete conformal structure in $\mathcal{U}^H$ with combinatorial $\alpha$-curvature $\overline{R}$.
  \item[(2)] Suppose that there exists a discrete conformal structure $\overline{u}\in \mathcal{U}^H$ with combinatorial $\alpha$-curvature $\overline{R}$ and $\alpha \overline{R}\leq 0$,  then there exists a constant $\delta >0$ such that if $||R_\alpha(u(0))-R_\alpha(\overline{u})||<\delta$, the solutions of modified combinatorial $\alpha$-Ricci flow (\ref{modified combinatorial alpha Ricci flow}) and modified combinatorial $\alpha$-Calabi flow (\ref{modified combinatorial alpha Calabi flow}) exist for all time and converge exponentially fast to $\overline{u}$ respectively.
\end{description}
\end{theorem}
\proof
The proof of Theorem \ref{alpha flows with small initial energy 2} is similar to that of Theorem \ref{alpha flows with small initial energy 1}, so we only give some necessary steps of the proof for the second part here.
For the modified combinatorial $\alpha$-Ricci flow (\ref{modified combinatorial alpha Ricci flow}), set $\Gamma_i(u)=\overline{R}_i-R_{\alpha,i}$.
By direct calculations, we have
\begin{equation*}
D\Gamma|_{u=\overline{u}}=-\Sigma^{-\alpha}\Lambda^H+\alpha L' =-\Sigma^{-\frac{\alpha}{2}}
(\Sigma^{-\frac{\alpha}{2}}\Lambda^H\Sigma^{-\frac{\alpha}{2}}-\alpha L')\Sigma^{\frac{\alpha}{2}},
\end{equation*}
where $\Sigma=\text{diag}\{e^{u_1},e^{u_2},...,e^{u_N}\}$, $L'=\text{diag}\{\overline{R}_1,\overline{R}_2,...,\overline{R}_N\}$.
By Lemma \ref{hyperbolic extended lemma} and the condition $\alpha \overline{R}\leq 0$, $D\Gamma|_{u=\overline{u}}$ has $N$ negative eigenvalues, which implies $\overline{u}$ is a local attractor of the modified combinatorial $\alpha$-Ricci flow (\ref{modified combinatorial alpha Ricci flow}). Then the conclusion follows from Lyapunov Stability Theorem (\cite{Pontryagin}, Chapter 5).

Similarly, for the modified combinatorial $\alpha$-Calabi flow (\ref{modified combinatorial alpha Calabi flow}), set $\Gamma_i(u)=\Delta(R_\alpha-\overline{R})_i$.
By direct calculations, we have
\begin{equation*}
\begin{aligned}
D\Gamma|_{u=\overline{u}}
&=-\Sigma^{-\alpha}\Lambda^H\Sigma^{-\alpha}\Lambda^H+\alpha\Sigma^{-\alpha}\Lambda^H L' \\
&=-\Sigma^{-\frac{\alpha}{2}}\left(\Sigma^{-\frac{\alpha}{2}}\Lambda^H\Sigma^{-\alpha}\Lambda^H\Sigma^{-\frac{\alpha}{2}}
-\alpha\Sigma^{-\frac{\alpha}{2}} \Lambda^H\Sigma^{-\frac{\alpha}{2}} L' \right)\Sigma^{\frac{\alpha}{2}}\\
&=-\Sigma^{-\frac{\alpha}{2}}\left(Q^2-\alpha Q L'\right)\Sigma^{\frac{\alpha}{2}}\\
&=-\Sigma^{-\frac{\alpha}{2}}Q^{\frac{1}{2}}\left(Q^2+Q^{\frac{1}{2}}(-\alpha L')Q^{\frac{1}{2}}\right)Q^{-\frac{1}{2}}\Sigma^{\frac{\alpha}{2}},
\end{aligned}
\end{equation*}
where $Q=\Sigma^{-\frac{\alpha}{2}} \Lambda^H\Sigma^{-\frac{\alpha}{2}}$ is a symmetric and positive definite matrix by Lemma \ref{hyperbolic extended lemma}. By $\alpha \overline{R}\leq 0$, $D\Gamma|_{u=\overline{u}}$ has $N$ negative eigenvalues, which implies $\overline{u}$ is a local attractor of the modified combinatorial $\alpha$-Calabi flow (\ref{modified combinatorial alpha Calabi flow}). Then the conclusion follows from Lyapunov Stability Theorem (\cite{Pontryagin}, Chapter 5).
\qed

\subsection{Uniqueness for the solution of extended combinatorial $\alpha$-Ricci flow}
Theorem \ref{alpha flows with small initial energy 1} and Theorem \ref{alpha flows with small initial energy 2} give the longtime existence and convergence for the solutions of the combinatorial $\alpha$-curvature flows for initial values with small energy.
But for general initial value, the combinatorial $\alpha$-curvature flows may develop singularities, including the removable singularities and the essential singularities.
For the combinatorial $\alpha$-Ricci flow, one can extend it through the removable singularities by extending the combinatorial $\alpha$-curvature.

\begin{definition}\label{extended combinatorial alpha Ricci flow definition}
Suppose $(M,\mathcal{T},\varepsilon,\eta)$ is a weighted triangulated connected closed surface with the weights $\varepsilon: V\rightarrow \{0,1\}$ and $\eta: E\rightarrow \mathbb{R}$ satisfying the structure conditions (\ref{structure condition 1}) and (\ref{structure condition 2}). For Euclidean and hyperbolic background geometry, the extended combinatorial $\alpha$-Ricci flow for discrete conformal structures on polyhedral surfaces is defined to be
\begin{equation*}
\frac{du_i}{dt}=-\widetilde{R}_{\alpha,i}
\end{equation*}
and the extended modified combinatorial $\alpha$-Ricci flow is defined to be
\begin{equation}\label{extended modified combinatorial alpha Ricci flow}
\frac{du_i}{dt}=\overline{R}_i-\widetilde{R}_{\alpha,i},
\end{equation}
where $\overline{R}: V\rightarrow \mathbb{R}$ is a function and $\widetilde{R}_{\alpha,i}=\frac{\widetilde{K}_i}{e^{\alpha u_i}}=\frac{1}{e^{\alpha u_i}}(2\pi-\sum_{\{ijk\}\in F}\widetilde{\theta}_i^{jk})$.
\end{definition}

Note that for the modified combinatorial $\alpha$-Ricci flow (\ref{modified combinatorial alpha Ricci flow}),
$\overline{R}_i-R_{\alpha,i}$  is smooth and locally Lipschitz as a function of $u\in \mathcal{U}$ ($\mathcal{U}^E$ or $\mathcal{U}^H$).
By Picard's uniqueness for the solution of ODE, the modified combinatorial $\alpha$-Ricci flow (\ref{modified combinatorial alpha Ricci flow}) has a unique solution $u(t)$, $t\in[0,T)$ for some $T>0$.
As noted by the first author in \cite{Xu 21}, the extension $\widetilde{K}$ of the classical combinatorial curvature $K$ is not a locally
 Lipschitz function of the generalized discrete conformal structures.
 This implies that the extended combinatorial $\alpha$-curvature $\widetilde{R}_{\alpha}$ is only a continuous function of the generalized discrete conformal structures $f\in \mathbb{R}^N$ and not a locally Lipschitz function. Therefore, there may exist more than one solution
  for the extended combinatorial $\alpha$-Ricci flow by the standard ODE theory. However, we can prove the following uniqueness for the solutions of the extended modified combinatorial $\alpha$-Ricci flow (\ref{extended modified combinatorial alpha Ricci flow})
  under the structure conditions (\ref{structure condition 1}) and (\ref{structure condition 2}), which is a generalization of Theorem \ref{main theorem 1} (a).

\begin{theorem}\label{unique theorem}
Suppose $(M,\mathcal{T},\varepsilon,\eta)$ is a weighted triangulated connected closed surface with the weights $\varepsilon: V\rightarrow \{0,1\}$ and $\eta: E\rightarrow \mathbb{R}$ satisfying the structure conditions (\ref{structure condition 1}) and (\ref{structure condition 2}). $\alpha\in \mathbb{R}$ is a constant and $\overline{R}: V\rightarrow \mathbb{R}$ is a given function with $\alpha\overline{R}\leq 0$. For any initial generalized discrete conformal structure, the solution of the extended modified combinatorial $\alpha$-Ricci flow (\ref{extended modified combinatorial alpha Ricci flow}) is unique.
\end{theorem}
\proof
In the case of $\alpha=0$, Theorem \ref{unique theorem} was proved by the first author in \cite{Xu 21}.
Here we only prove the case of $\alpha\neq0$, the idea of which comes from \cite{Ge-Hua 20}.
Define the following map in Euclidean background geometry
\begin{equation*}
\begin{aligned}
w^E(u) :\  &\mathbb{R}^N\rightarrow W^E\subseteq\mathbb{R}^N\\
&u\mapsto w^E(u):=(w^E_1(u_1),...,w^E_N(u_N)),
\end{aligned}
\end{equation*}
where $w^E_i(u_i)=\int^{u_i}e^{\frac{1}{2}\alpha x}dx$ and $W^E= w^E(\mathbb{R}^N)$.
Note that $w^E_i(u_i)$ is strictly increasing in $u_i$, the map $w^E(u): \mathbb{R}^N\rightarrow W^E$ is a diffeomorphism.
The map $w^H(u)$ in hyperbolic background geometry can be defined similarly. Except for different domains and images, the map $w^E_i$ and $w^H_i$ are similar, we use $w_i$ to represent $w^E_i$ and $w^H_i$ for simplicity, if there is no confusion. Denote the inverse map of $w(u)$ by $u(w)$, we define a new function on $W$ ($W^E$ or $W^H$) by
\begin{equation*}
\widehat{F}(w):=\widetilde{F}(u(w)),
\end{equation*}
where $\widetilde{F}(u)$ is defined by (\ref{Euclidean extended F}) for Euclidean background geometry and by (\ref{hyperbolic extended F}) for hyperbolic background geometry. Note that the extended modified combinatorial $\alpha$-Ricci flow (\ref{extended modified combinatorial alpha Ricci flow}) can be written as
\begin{equation}\label{negative gradient flow}
u'_i(t)=-\frac{1}{e^{\alpha u_i}}\nabla_{u_i}\widetilde{F}.
\end{equation}
Therefore, the equation (\ref{negative gradient flow}) is equivalent to the following equation in the $w$-coordinate
\begin{equation*}
w'_i(t)=-\nabla_{w_i}\widehat{F},
\end{equation*}
which means that the extended modified combinatorial $\alpha$-Ricci flow (\ref{extended modified combinatorial alpha Ricci flow}) is equivalent to a negative gradient flow of the function $\widehat{F}$ in the $w$-coordinate.

Suppose $u_A(t)$ and $u_B(t)$ are two solutions of the extended modified combinatorial $\alpha$-Ricci flow (\ref{extended modified combinatorial alpha Ricci flow}) with $u_A(0)=u_B(0)$.
Then $w_A(0)=w_B(0)$ as the map $w(u)$ is a diffeomorphism. There exists $T\in(0,+\infty)$ such that
$\{u_A(t)|t\in[0,T]\}$ and $\{u_B(t)|t\in[0,T]\}$ lie in a compact set and hence $\{w_A(t)|t\in[0,T]\}$ and $\{w_B(t)|t\in[0,T]\}$ lie in a compact convex subsect $W'$ of $W$. We claim that $\widehat{F}$ is $(-\lambda)$-convex on $W'$, i.e. there exists a positive constant $\lambda$ such that for any $w_A, w_B\in W'$,
\begin{equation*}
(\nabla_w\widehat{F}(w_A)-\nabla_w\widehat{F}(w_B))\cdot (w_A-w_B)+\lambda |w_A-w_B|^2\geq0.
\end{equation*}
Set $h(t)=|w_A(t)-w_B(t)|^2$, then
\begin{equation*}
\begin{aligned}
h'(t)
&=2(w_A(t)-w_B(t))\cdot (w'_A(t)-w'_B(t))\\
&=-2(w_A(t)-w_B(t))\cdot (\nabla_w\widehat{F}(w_A(t))-\nabla_w\widehat{F}(w_B(t)))\\
&\leq2\lambda |w_A(t)-w_B(t)|^2=2\lambda h(t)
\end{aligned}
\end{equation*}
by the claim, which implies $h(t)\leq h(0)e^{2\lambda t}$. Note that $h(0)=0$ and $h(t)\geq0$, then we have $h(t)\equiv0$.
Therefore, $w_A(t)=w_B(t)$, which implies $u_A(t)=u_B(t)$.

We use Ge-Hua's trick in \cite{Ge-Hua 20} to prove the claim. Note that $\widetilde{F}$ is a $C^1$-smooth convex function defined on $\mathbb{R}^N$ by the proof of  Theorem \ref{Euclidean global rigidity theorem} for Euclidean background geometry and by the proof of  Theorem \ref{hyperbolic global rigidity theorem} for hyperbolic background geometry.
By mollifying $\widetilde{F}$ using the standard mollifier $\varphi_\epsilon(u)=\frac{1}{\epsilon^N}\varphi(\frac{u}{\epsilon})$ with
\begin{eqnarray*}
\varphi(u)=
\begin{cases}
Ce^{\frac{1}{1-|u|^2}},\ &|u|<1, \\
0,\ &|u|\geq1,
\end{cases}
\end{eqnarray*}
where the positive constant $C$ is chosen such that $\int_{\mathbb{R}^N}\varphi(u)du=1$, we have $\widetilde{F}_\epsilon=\widetilde{F}\ast \varphi_\epsilon$ is a smooth convex function of $u$ and $\widetilde{F}_\epsilon\rightarrow \widetilde{F}$ in $C^1_{loc}$ as $\epsilon\rightarrow 0$.
Moreover, $\nabla_{u_i} \widetilde{F}_\epsilon=\nabla_{u_i} \widetilde{F}\ast \varphi_\epsilon=(\widetilde{K}_i-\overline{R}_ie^{\alpha u_i})\ast \varphi_\epsilon.$
Set $\widehat{F}_\epsilon(w)=\widetilde{F}_\epsilon(u(w))$, then by the chain rules, we have $\nabla_w\widehat{F}_\epsilon=\nabla_u\widetilde{F}_\epsilon\frac{\partial u}{\partial w}$ and
\begin{equation*}
\begin{aligned}
\nabla^2_{w_iw_j}\widehat{F}_\epsilon
&=\frac{\partial^2\widetilde{F}_\epsilon}{\partial u_i\partial u_j}\frac{\partial u_i}{\partial w_i}\frac{\partial u_j}{\partial w_j}+\frac{\partial \widetilde{F}_\epsilon}{\partial u_i}\frac{\partial^2u_i}{\partial w_i\partial w_j}\\
&=\frac{\partial^2\widetilde{F}_\epsilon}{\partial u_i\partial u_j}\frac{\partial u_i}{\partial w_i}\frac{\partial u_j}{\partial w_j}-\frac{1}{2}(\widetilde{K}_i-\overline{R}_ie^{\alpha u_i})\ast \varphi_\epsilon\alpha e^{-\alpha u_i}\delta_{ij}.
\end{aligned}
\end{equation*}
Note that on the compact subset $W'$, there exists a positive constant $\lambda$ such that
\begin{equation*}
\frac{1}{2}|(\widetilde{K}_i-\overline{R}_ie^{\alpha u_i})\ast \varphi_\epsilon\alpha e^{-\alpha u_i}|\leq C|\alpha|e^{-\alpha u_i}\leq\lambda.
\end{equation*}
By the convexity of $\widetilde{F}_\epsilon$, we have
\begin{equation*}
\nabla^2_w \widehat{F}_\epsilon\geq\frac{\partial u}{\partial w}\nabla^2_u\widetilde{F}_\epsilon\frac{\partial u}{\partial w}-\lambda I\geq-\lambda I,
\end{equation*}
which implies that $\widehat{F}_\epsilon$ is $(-\lambda)$-convex on $W'$, i.e. for any $w_A, w_B\in W'$,
\begin{equation*}
(\nabla_w\widehat{F}_\epsilon(w_A)-\nabla_w\widehat{F}_\epsilon(w_B))\cdot (w_A-w_B)+\lambda |w_A-w_B|^2\geq0.
\end{equation*}
The claim follows by letting $\epsilon\rightarrow0$.
\qed

\begin{remark}
By Definition \ref{extended combinatorial alpha Ricci flow definition}, $\widetilde{R}_\alpha|_{\mathcal{U}}=R_\alpha$ and the solutions of (\ref{modified combinatorial alpha Ricci flow}) and (\ref{extended modified combinatorial alpha Ricci flow}) agree on $\mathcal{U}$ by Picard's uniqueness for the solution of ODE. Therefore, any solution of the extended modified combinatorial $\alpha$-Ricci flow (\ref{extended modified combinatorial alpha Ricci flow}) extends the solution of the modified combinatorial $\alpha$-Ricci flow (\ref{modified combinatorial alpha Ricci flow}).
\end{remark}

\begin{remark}
By the proof of Theorem \ref{unique theorem}, one can further prove the uniqueness of the solution of extended modified combinatorial $\alpha$-Ricci flow (\ref{extended modified combinatorial alpha Ricci flow}) without the assumption $\alpha\overline{R}\leq 0$.
\end{remark}

\subsection{Longtime existence and global convergence of the combinatorial $\alpha$-Ricci flow in the Euclidean background geometry}
According to Definition \ref{extended combinatorial alpha Ricci flow definition}, the extended normalized combinatorial $\alpha$-Ricci flow for discrete conformal structures in the Euclidean background geometry is defined to be
\begin{equation}\label{extended normalized combinatorial alpha Ricci flow}
\frac{du_i}{dt}=R_{\alpha,av}-\widetilde{R}_{\alpha,i},
\end{equation}
where $R_{\alpha,av}=\frac{2\pi \chi(M)}{\sum^N_{i=1}e^{\alpha u_i}}$ is the average combinatorial $\alpha$-curvature.
We have the following result on the longtime existence and convergence for the solution of the normalized combinatorial $\alpha$-Ricci flow under the existence of the Euclidean discrete conformal structure with constant combinatorial $\alpha$-curvature, which is a slight generalization of Theorem \ref{main theorem 1} (b).

\begin{theorem}\label{Euclidean exist for all time and converge}
Suppose $(M,\mathcal{T},\varepsilon,\eta)$ is a weighted triangulated connected closed surface with the weights $\varepsilon: V\rightarrow \{0,1\}$ and $\eta: E\rightarrow \mathbb{R}$ satisfying the structure conditions (\ref{structure condition 1}) and (\ref{structure condition 2}).
$\alpha\in \mathbb{R}$ is a constant and $\alpha\chi(M)\leq 0$.
Suppose there exists a Euclidean discrete conformal structure $\overline{u}\in \mathcal{U}^E$ with constant combinatorial $\alpha$-curvature.
Then the normalized combinatorial $\alpha$-Ricci flow (\ref{normalized combinatorial alpha Ricci flow}) develops no essential singularities.
If the solution of (\ref{normalized combinatorial alpha Ricci flow}) develops no removable singularities in finite time, then the solution of (\ref{normalized combinatorial alpha Ricci flow}) exists for all time, converges exponentially fast to $\overline{u}$ for any initial Euclidean discrete conformal structure $u(0)\in  \mathcal{U}^E$ with
$\sum^N_{i=1}e^{\alpha u(0)}=\sum^N_{i=1}e^{\alpha\overline{u}_i}$
for $\alpha\neq 0$ and
$\sum^N_{i=1}u(0)=\sum^N_{i=1}\overline{u}_i$
for $\alpha=0$, and does not develop removable singularities at time infinity.
Furthermore, the unique solution of the extended normalized combinatorial $\alpha$-Ricci flow (\ref{extended normalized combinatorial alpha Ricci flow}) exists for all time and converges exponentially fast to $\overline{u}$ for any initial generalized Euclidean discrete conformal structure
$u(0)\in \mathbb{R}^N$ with $\sum^N_{i=1}e^{\alpha u(0)}=\sum^N_{i=1}e^{\alpha\overline{u}_i}$ for $\alpha\neq 0$ and $\sum^N_{i=1}u(0)=\sum^N_{i=1}\overline{u}_i$ for $\alpha=0$.
\end{theorem}
\proof
If $\alpha=0$, Theorem \ref{Euclidean exist for all time and converge} has been proved by the first author in \cite{Xu 21}.
Here we only prove the case of $\alpha\neq0$.
Similar to the proof of Lemma \ref{invariant along flow}, $\sum^N_{i=1}r^\alpha_i$ is invariant along the extended normalized combinatorial $\alpha$-Ricci flow (\ref{extended normalized combinatorial alpha Ricci flow}). Without loss of generality, assume $\sum^N_{i=1}r^\alpha_i(0)=N$, then the solution of the extended normalized combinatorial $\alpha$-Ricci flow (\ref{extended normalized combinatorial alpha Ricci flow}) stays in the hypersurface $\Pi:=\{u\in \mathbb{R}^N|\sum^N_{i=1}r^\alpha_i(0)=N\}$.

Suppose there exists a Euclidean discrete conformal structure $\overline{u}\in \mathcal{U}^E\cap \Pi$ with constant combinatorial $\alpha$-curvature $R_{\alpha,av}$. Define the following Ricci energy function
\begin{equation*}
\widetilde{F}(u)=-\sum_{\{ijk\}\in F}\widetilde{F}_{ijk}(u_i,u_j,u_k)
+\int^u_{\overline{u}}\sum_{i=1}^{N}(2\pi-R_{\alpha,av}r^{\alpha}_i)du_i.
\end{equation*}
By direct calculations, we have
\begin{equation*}
\nabla_{u_i}\widetilde{F}(u)=-\sum_{\{ijk\}\in F}\widetilde{\theta}_i^{jk}+2\pi-R_{\alpha,av}r^{\alpha}_i
=\widetilde{K}_i-R_{\alpha,av}r^{\alpha}_i
\end{equation*}
on $\mathbb{R}^N$ and
\begin{equation*}
\begin{aligned}
\mathrm{Hess}_u\widetilde{F}(u)
&=\Lambda^E-\alpha R_{\alpha,av}\Sigma^{\frac{\alpha}{2}}\left(I-\frac{r^{\frac{\alpha}{2}}\cdot (r^{\frac{\alpha}{2}})^T}{||r||^\alpha_\alpha}]\right)\Sigma^{\frac{\alpha}{2}}\\
&=\Sigma^{\frac{\alpha}{2}}\left(L-\alpha R_{\alpha,av}\frac{r^{\frac{\alpha}{2}}\cdot (r^{\frac{\alpha}{2}})^T}{||r||^\alpha_\alpha}]\right)\Sigma^{\frac{\alpha}{2}}
\end{aligned}
\end{equation*}
on $\mathcal{U}^E$, where $\Sigma=\text{diag}\{r_1,r_2,...,r_N\}$, $L=\Sigma^{-\frac{\alpha}{2}}\Lambda^E\Sigma^{-\frac{\alpha}{2}}$.
Note that $\widetilde{F}(u)$ is a $C^1$ smooth convex function defined on $\mathbb{R}^N$
and $\nabla_{u_i} \widetilde{F}(\overline{u})=\widetilde{K}_i(\overline{u})-R_{\alpha,av}e^{\alpha \overline{u}_i}=\widetilde{R}_{\alpha,i}(\overline{u})e^{\alpha \overline{u}_i}-R_{\alpha,av}e^{\alpha \overline{u}_i}=0$ by assumption,
we have $\widetilde{F}(u)\geq\widetilde{F}(\overline{u})=0$ and $\lim_{u\rightarrow \infty}\widetilde{F}(u)|_\Pi=+\infty$ by the following property of convex functions, a proof of which could be found in \cite{Ge-Xu 18} (Lemma 4.6).
\begin{lemma}\label{Xu Lemma 4.6}
Suppose $f(x)$ is a $C^1$ smooth convex function on $\mathbb{R}^N$ with $\nabla f(x)=0$ for some $x_0\in \mathbb{R}^N$, $f(x)$ is $C^2$ smooth and strictly convex in a neighborhood of $x_0$, then $\lim_{x\rightarrow \infty}f(x)=+\infty$.
\end{lemma}
Note that
\begin{equation*}
\frac{d}{dt}\widetilde{F}(u(t))=\sum_{i=1}^N\nabla_{u_i}\widetilde{F}\cdot\frac{du_i}{dt}
=\sum_{i=1}^N(\widetilde{K}_i-R_{\alpha,av}r_i^\alpha)(R_{\alpha,av}-\widetilde{R}_{\alpha,i})
=-\sum_{i=1}^N(R_{\alpha,av}-\widetilde{R}_{\alpha,i})^2r_i^\alpha\leq0
\end{equation*}
along the extended normalized combinatorial $\alpha$-Ricci flow (\ref{extended normalized combinatorial alpha Ricci flow}), we have
$\widetilde{F}(u(t))$ is bounded along (\ref{extended normalized combinatorial alpha Ricci flow}), which implies
the solution $u(t)$ of the extended normalized combinatorial $\alpha$-Ricci flow (\ref{extended normalized combinatorial alpha Ricci flow}) stays in a compact subset of $\Pi$ by $\lim_{u\rightarrow \infty}\widetilde{F}(u)|_\Pi=+\infty$.
Therefore, the normalized combinatorial $\alpha$-Ricci flow (\ref{normalized combinatorial alpha Ricci flow}) develops no essential singularities and the solution of extended normalized combinatorial $\alpha$-Ricci flow (\ref{extended normalized combinatorial alpha Ricci flow}) exists for all time.
Furthermore, if the solution of (\ref{normalized combinatorial alpha Ricci flow}) develops no removable singularities in finite time,
$\lim_{t\rightarrow +\infty}\widetilde{F}(u(t))$ exists and there exists $\xi_n\in (n,n+1)$ such that
\begin{equation}\label{1}
\widetilde{F}(u(n+1))-\widetilde{F}(u(n))=\frac{d}{dt}\widetilde{F}(u(t))|_{t=\xi_n}
=-\sum_{i=1}^N(R_{\alpha,av}-\widetilde{R}_{\alpha,i})^2r_i^\alpha|_{t=\xi_n}\rightarrow 0,\ \ n\rightarrow \infty.
\end{equation}
As $u(\xi_n)$ is bounded, there exists a convergent subsequence of $u(\xi_n)$, still denoted by $u(\xi_n)$ for simplicity, such that $u(\xi_n)\rightarrow u^*\in \mathbb{R}^N$ and then $\widetilde{R}_\alpha(u(\xi_n))\rightarrow \widetilde{R}_\alpha(u^*)$.
The equation (\ref{1}) shows that $\widetilde{R}_\alpha(u^*)=R_{\alpha,av}=R_\alpha(\overline{u})$.
Therefore, $u^*=\overline{u}$ by Theorem \ref{Euclidean global rigidity theorem}.

Similar to Theorem \ref{alpha flows with small initial energy 1}, set $\Gamma_i(u)=R_{\alpha,av}-\widetilde{R}_{\alpha,i}$, then $D\Gamma|_{u=\overline{u}}$ restricted to the hypersurface $\Pi$ is negative definite,
which implies that $\overline{u}$ is a local attractor of (\ref{extended normalized combinatorial alpha Ricci flow}). Then the conclusion follows from Lyapunov Stability Theorem (\cite{Pontryagin}, Chapter 5).
\qed

For the prescribed combinatorial $\alpha$-curvature problem in the Euclidean background geometry, we have following result paralleling to Theorem \ref{Euclidean exist for all time and converge}.

\begin{theorem}\label{new theorem}
Suppose $(M,\mathcal{T},\varepsilon,\eta)$ is a weighted triangulated connected closed surface with the weights $\varepsilon: V\rightarrow \{0,1\}$ and $\eta: E\rightarrow \mathbb{R}$ satisfying the structure conditions (\ref{structure condition 1}) and (\ref{structure condition 2}).
$\alpha\in \mathbb{R}$ is a constant and $\overline{R}: V\rightarrow \mathbb{R}$ is a given function.
\begin{description}
  \item[(1)] If the solution of modified combinatorial $\alpha$-Ricci flow (\ref{modified combinatorial alpha Ricci flow}) in the Euclidean background geometry converges, then there exists a Euclidean discrete conformal structure $\overline{u}\in \mathcal{U}^E$ with combinatorial $\alpha$-curvature $\overline{R}$.
  \item[(2)] Suppose there exists a Euclidean discrete conformal structure $\overline{u}\in \mathcal{U}^E$ with combinatorial $\alpha$-curvature $\overline{R}$ satisfying $\alpha\overline{R}\leq 0$ and $\alpha\overline{R}\not\equiv 0$.
      Then the modified combinatorial $\alpha$-Ricci flow (\ref{modified combinatorial alpha Ricci flow}) in the Euclidean background geometry
      develops no essential singularities along the flow.
      If the solution of (\ref{modified combinatorial alpha Ricci flow}) develops no removable singularities in finite time,
      then the solution of (\ref{modified combinatorial alpha Ricci flow}) exists for all time,
      converges exponentially fast to $\overline{u}$ for any initial Euclidean discrete conformal structure $u(0)\in \mathcal{U}^E$
      and does not develop removable singularity at time infinity.
      Furthermore, the unique solution of the extended modified combinatorial $\alpha$-Ricci flow (\ref{extended modified combinatorial alpha Ricci flow}) exists for all time and converges exponentially fast to $\overline{u}$ for any initial generalized Euclidean discrete conformal structure $u(0)\in \mathbb{R}^N$.
\end{description}
\end{theorem}
\proof
The proof of (1) in Theorem \ref{new theorem} is similar to that of (1) in Theorem \ref{alpha flows with small initial energy 1}, so we omit the details of the proof here.
By Theorem \ref{Euclidean global rigidity theorem}, the function
\begin{equation*}
\widetilde{F}(u)=-\sum_{\{ijk\}\in F}\widetilde{F}_{ijk}(u_i,u_j,u_k)
+\int^u_{\overline{u}}\sum_{i=1}^{N}
(2\pi-\overline{R}_ir^{\alpha}_i)du_i
\end{equation*}
is a $C^1$-smooth convex function defined on $\mathbb{R}^N$ and locally strictly convex on $\mathcal{U}^E$ in the case of $\alpha\overline{R}\leq 0$ and $\alpha\overline{R}\not\equiv 0$.
By direct calculations, we have
\begin{equation*}
\nabla_{u_i}\widetilde{F}(u)=-\sum_{\{ijk\}\in F}\widetilde{\theta}_i^{jk}+2\pi-\overline{R}_ir^{\alpha}_i
=\widetilde{K}_i-\overline{R}_ir^{\alpha}_i.
\end{equation*}
By assumption,
$\nabla_{u_i}\widetilde{F}(\overline{u})=\widetilde{K}_i(\overline{u})-\overline{R}_ie^{\alpha \overline{u}_i}=\widetilde{R}_{\alpha,i}(\overline{u})e^{\alpha \overline{u}_i}-\overline{R}_ie^{\alpha \overline{u}_i}=0.$
By Lemma \ref{Xu Lemma 4.6}, since $\widetilde{F}(u)$ is strictly convex in a neighborhood of $\overline{u}\in \mathcal{U}^E$, we have $\lim_{u\rightarrow \infty}\widetilde{F}(u)=+\infty$.
Note that
\begin{equation*}
\frac{d}{dt}\widetilde{F}(u(t))=\sum_{i=1}^N\nabla_{u_i}\widetilde{F}\cdot\frac{du_i}{dt}
=\sum_{i=1}^N(\widetilde{K}_i-\overline{R}_ir_i^\alpha)(\overline{R}_i-\widetilde{R}_{\alpha,i})
=-\sum_{i=1}^N(\overline{R}_i-\widetilde{R}_{\alpha,i})^2r_i^\alpha\leq0
\end{equation*}
along the extended modified combinatorial $\alpha$-Ricci flow (\ref{extended modified combinatorial alpha Ricci flow}), we have
$\widetilde{F}(u(t))$ is bounded,
which implies that the solution $u(t)$ of the extended modified combinatorial $\alpha$-Ricci flow (\ref{extended modified combinatorial alpha Ricci flow}) stays in a compact subset of $\mathbb{R}^N$ and $\lim_{t\rightarrow +\infty}\widetilde{F}(u(t))$ exists.
Therefore, the solution of the extended modified combinatorial $\alpha$-Ricci flow (\ref{extended modified combinatorial alpha Ricci flow})
exists for all time and the modified combinatorial $\alpha$-Ricci flow (\ref{modified combinatorial alpha Ricci flow}) in the Euclidean background geometry develops no essential singularities.
By the mean value theorem, there exists $\xi_n\in (n,n+1)$ such that
\begin{equation}\label{3}
\widetilde{F}(u(n+1))-\widetilde{F}(u(n))=\frac{d}{dt}\widetilde{F}(u(t))|_{t=\xi_n}
=-\sum_{i=1}^N(\overline{R}_i-\widetilde{R}_{\alpha,i})^2r_i^\alpha|_{t=\xi_n}\rightarrow 0,\ \ n\rightarrow \infty.
\end{equation}
As $\{u(\xi_n)\}$ is bounded, there exists a convergent subsequence of $\{u(\xi_n)\}$, still denoted by $\{u(\xi_n)\}$ for convenience, such that $u(\xi_n)\rightarrow u^*\in \mathbb{R}^N$ and then $\widetilde{R}_\alpha(u(\xi_n))\rightarrow \widetilde{R}_\alpha(u^*)$.
The equation (\ref{3}) shows that $\widetilde{R}_\alpha(u^*)=\overline{R}=R_\alpha(\overline{u})$.
Therefore, $u^*=\overline{u}$ by Theorem \ref{Euclidean global rigidity theorem}.

Set $\Gamma_i(u)=\overline{R}_i-\widetilde{R}_{\alpha,i}$, then
\begin{equation*}
D\Gamma|_{u=\overline{u}}=-\Sigma^{-\alpha}\Lambda^E+\alpha L' =-\Sigma^{-\frac{\alpha}{2}}
(\Sigma^{-\frac{\alpha}{2}}\Lambda^E\Sigma^{-\frac{\alpha}{2}}-\alpha L')\Sigma^{\frac{\alpha}{2}},
\end{equation*}
where $\Sigma=\text{diag}\{e^{u_1},e^{u_2},...,e^{u_N}\}$, $L'=\text{diag}\{\overline{R}_1,\overline{R}_2,...,\overline{R}_N\}$. Therefore,
$D\Gamma|_{u=\overline{u}}$ is negative definite by $\alpha\overline{R}\leq 0$ and $\alpha\overline{R}\not\equiv 0$,
which implies that $\overline{u}$ is a local attractor of (\ref{extended modified combinatorial alpha Ricci flow}). Then the conclusion follows from Lyapunov Stability Theorem (\cite{Pontryagin}, Chapter 5).
\qed

\subsection{Longtime existence and global convergence of the combinatorial $\alpha$-Ricci flow in the hyperbolic background geometry}
In the hyperbolic background geometry, we take the constant combinatorial $\alpha$-curvature problem as a special case of the prescribed combinatorial $\alpha$-curvature problem. To study the prescribed combinatorial $\alpha$-curvature problem in the hyperbolic background geometry using the extended modified combinatorial $\alpha$-Ricci flow (\ref{extended modified combinatorial alpha Ricci flow}),
we need the following lemma, which was proved by the first author in \cite{Xu 21} (Corollary 4.9).

\begin{lemma}(\cite{Xu 21})\label{Xu Corollary 4.9}
Suppose $(M,\mathcal{T},\varepsilon,\eta)$ is a weighted triangulated connected closed surface with the weights $\varepsilon: V\rightarrow \{0,1\}$ and $\eta: E\rightarrow \mathbb{R}$ satisfying the structure conditions (\ref{structure condition 1}) and (\ref{structure condition 2}).
$\{ijk\}\in F$ is a triangle with $\varepsilon_i=1$. Then for any $\epsilon>0$, there exists a positive number $L=L(\varepsilon, \eta, \epsilon)$ such that if $f_i>L$, the extended inner $\widetilde{\theta}_i^{jk}$ at the vertex $i\in V$ in the generalized hyperbolic triangle $\{ijk\}\in F$ with edge lengths given by (\ref{hyperbolic edge length}) is smaller than $\epsilon$.
\end{lemma}

We have the following result on the longtime existence and convergence of the solution of extended modified combinatorial $\alpha$-Ricci flow (\ref{extended modified combinatorial alpha Ricci flow}) in hyperbolic background geometry, which is a generalization of
Theorem \ref{main theorem 1} (c).

\begin{theorem}\label{hyperbolic exist for all time}
Suppose $(M,\mathcal{T},\varepsilon,\eta)$ is a weighted triangulated connected closed surface with the weights $\varepsilon: V\rightarrow \{0,1\}$ and $\eta: E\rightarrow \mathbb{R}$ satisfying the structure conditions (\ref{structure condition 1}) and (\ref{structure condition 2}). $\alpha\in \mathbb{R}$ is a constant and $\overline{R}: V\rightarrow \mathbb{R}$ is a given function.
Suppose there exists a hyperbolic discrete conformal structure $\overline{u}\in \mathcal{U}^H$ with combinatorial $\alpha$-curvature $\overline{R}$ such that one of the following three conditions is satisfied
\begin{description}
  \item[(1)] $\alpha>0$ and $\overline{R}_i\leq0$ for all $i\in V$,
  \item[(2)] $\alpha<0$ and $\overline{R}_i\in [0,2\pi)$ for all $i\in V$,
  \item[(3)] (\cite{Xu 21}) $\alpha=0$, $\overline{R}_i\in (-\infty, 2\pi)$ for all $i\in V$ and $\sum^N_{i=1} \overline{R}_{i}>2\pi \chi(M)$.
\end{description}
Then the modified combinatorial $\alpha$-Ricci flow (\ref{modified combinatorial alpha Ricci flow}) in the hyperbolic background geometry develops no essential singularities.
If the solution of (\ref{modified combinatorial alpha Ricci flow}) develops no removable singularities in finite time, then the solution of (\ref{modified combinatorial alpha Ricci flow}) exists for all time, converges exponentially fast to $\overline{u}$ for any initial hyperbolic  discrete conformal structure $u(0)\in \mathcal{U}^H$ and does not develop removable singularities at time infinity.
Furthermore, the unique solution of the extended modified combinatorial $\alpha$-Ricci flow (\ref{extended modified combinatorial alpha Ricci flow}) exists for all time and converges exponentially fast to $\overline{u}$ for any initial generalized hyperbolic discrete conformal structure $u(0)$.
\end{theorem}
\proof
The case of $\alpha=0$ has been proved by the first author in \cite{Xu 21}. As the proofs for the cases of $\alpha>0$ and $\alpha<0$ are all the same, we only prove the case of $\alpha>0$.

Suppose there exists a hyperbolic discrete conformal structure $\overline{u}\in \mathcal{U}^H$ with combinatorial $\alpha$-curvature $\overline{R}$, then the following Ricci energy function
\begin{equation*}
\widetilde{F}(u)=-\sum_{\{ijk\}\in F}\widetilde{F}_{ijk}(u_i,u_j,u_k)
+\int^u_{\overline{u}}\sum_{i=1}^{N}(2\pi-\overline{R}_ie^{\alpha u_i})du_i
\end{equation*}
is a $C^1$-smooth convex function defined on $\mathbb{R}^{N_0}\times \mathbb{R}^{N_1}_{<0}$ by the proof of Theorem \ref{hyperbolic global rigidity theorem}, where $N_0$ is the number of the vertices $i\in V$ with $\varepsilon_i=0$ and $N_1=N-N_0$. Moreover, $\nabla_{u_i} \widetilde{F}(\overline{u})=\widetilde{K}_i(\overline{u})-\overline{R}_ie^{\alpha \overline{u}_i}=\widetilde{R}_{\alpha,i}(\overline{u})e^{\alpha \overline{u}_i}-\overline{R}_ie^{\alpha \overline{u}_i}=0$ by assumption, which implies $\widetilde{F}(u)\geq\widetilde{F}(\overline{u})=0$ and $\lim_{u\rightarrow \infty}\widetilde{F}(u)=+\infty$ by Lemma \ref{Xu Lemma 4.6}.
Furthermore,
\begin{equation*}
\frac{d}{dt}\widetilde{F}(u(t))=\sum_{i=1}^N\nabla_{u_i}\widetilde{F}\cdot\frac{du_i}{dt}
=\sum_{i=1}^N(\widetilde{K}_i-\overline{R}_ie^{\alpha u_i})(\overline{R}_i-\widetilde{R}_{\alpha,i})
=-\sum_{i=1}^N(\overline{R}_i-\widetilde{R}_{\alpha,i})^2e^{\alpha u_i}\leq0
\end{equation*}
along the extended modified combinatorial $\alpha$-Ricci flow (\ref{extended modified combinatorial alpha Ricci flow}),
which implies the solution $u(t)$ of the extended modified combinatorial $\alpha$-Ricci flow (\ref{extended modified combinatorial alpha Ricci flow}) is bounded in $\mathbb{R}^N$. By Remark \ref{remark on u}, $u=(u_1,..,u_N)\in \mathbb{R}^{N_0}\times \mathbb{R}^{N_1}_{<0}$.
We claim that $u_i(t)$ is uniformly bounded from above in $\mathbb{R}_{<0}$ for vertices $i\in V$ with $\varepsilon_i=1$.
We shall prove the theorem assuming the claim and then prove the claim.

By the claim, we have the solution $u(t)$ of the extended modified combinatorial $\alpha$-Ricci flow (\ref{extended modified combinatorial alpha Ricci flow}) lies in a compact subset of $\mathbb{R}^{N_0}\times \mathbb{R}^{N_1}_{<0}$, which implies that the modified combinatorial $\alpha$-Ricci flow (\ref{modified combinatorial alpha Ricci flow}) in the hyperbolic background geometry develops no essential singularities and
the solution of the extended modified combinatorial $\alpha$-Ricci flow (\ref{extended modified combinatorial alpha Ricci flow}) exists for all time.
As $u(t)$ is bounded, $\widetilde{F}(u(t))$ is bounded along the extended modified combinatorial $\alpha$-Ricci flow (\ref{extended modified combinatorial alpha Ricci flow}) and $\lim_{t\rightarrow +\infty}\widetilde{F}(u(t))$ exists.
By the mean value theorem, there exists $\xi_n\in (n,n+1)$ such that
\begin{equation}\label{2}
\widetilde{F}(u(n+1))-\widetilde{F}(u(n))=\frac{d}{dt}\widetilde{F}(u(t))|_{t=\xi_n}
=-\sum_{i=1}^N(\overline{R}_i-\widetilde{R}_{\alpha,i})^2e^{\alpha u_i}|_{t=\xi_n}\rightarrow 0,\ \ n\rightarrow \infty.
\end{equation}
As $u(\xi_n)$ is bounded, then there exists a subsequence of $u(\xi_n)$, still denoted by $u(\xi_n)$ for simplicity,
such that $u(\xi_n)\rightarrow u^*\in \mathbb{R}^{N_0}\times \mathbb{R}^{N_1}_{<0}$ and then $\widetilde{R}_\alpha(u(\xi_n))\rightarrow \widetilde{R}_\alpha(u^*)$. The equation (\ref{2}) shows that $\widetilde{R}_\alpha(u^*)=\overline{R}=R_{\alpha}(\overline{u})$. Therefore, $u^*=\overline{u}$ by Theorem \ref{hyperbolic global rigidity theorem}.
Similar to Theorem \ref{alpha flows with small initial energy 2}, set $\Gamma_i(u)=\overline{R}_i-\widetilde{R}_{\alpha,i}$, then $D\Gamma|_{u=\overline{u}}$ is negative definite, which implies that $\overline{u}$ is a local attractor of (\ref{extended modified combinatorial alpha Ricci flow}). Then the conclusion follows from Lyapunov Stability Theorem (\cite{Pontryagin}, Chapter 5).

We use Ge-Xu's trick in \cite{Ge-Xu 17} to prove the claim. Suppose that
there exists at least one vertex $i\in V$ such that $\lim_{t\rightarrow T}u_i(t)=0$ for $T\in(0,+\infty)$, which is equivalent to $\lim_{t\rightarrow T}f_i(t)=+\infty$ by the map (\ref{hyperblic DCS in u}). For the vertex $i$, by Lemma \ref{Xu Corollary 4.9}, there exists $c\in (-\infty,0)$ such that whenever $u_i(t)>c$, the extended inner angle $\widetilde{\theta}_i$ is smaller than $\epsilon=\frac{2\pi}{d_i}>0$, where $d_i$ is the degree at the vertex $i$.
Thus $\widetilde{K}_i=2\pi-\sum_{\{ijk\}\in F}\widetilde{\theta}_i>0\geq\overline{R}_ie^{\alpha u_i}$ and $\widetilde{R}_{\alpha,i}=\frac{\widetilde{K}_i}{e^{\alpha u_i}}>\overline{R}_i$.
Choose a time $t_0\in (0,T)$ such that $u_i(t_0)>c$, this can be done because $\lim_{t\rightarrow T}u_i(t)=0$. Set $a=\inf\{t<t_0|u_i(s)>c, \forall s\in (t,t_0]\}$, then $u_i(a)=c$.
Note that $\frac{du_i}{dt}=\overline{R}_i-\widetilde{R}_{\alpha,i}<0$ on $(a,t_0]$, we have $u_i(t_0)<u_i(a)=c$, which contradicts $u_i(t_0)>c$.
Therefore, $u_i(t)$ is uniformly bounded from above in $\mathbb{R}_{<0}$ for vertex $i$ with $\varepsilon_i=1$.
\qed

\begin{remark}
If $\varepsilon_i=0$ for all $i\in V$, the condition (1)(2)(3) in Theorem \ref{hyperbolic exist for all time} can be generalized to be the condition that $\alpha\in \mathbb{R}$ is a constant and $\overline{R}: V\rightarrow \mathbb{R}$ is a given function defined on $V$ with $\alpha\overline{R}\leq 0$, which was proved in \cite{Xu Z}.
\end{remark}

\begin{remark}
By the proof of Theorem \ref{hyperbolic exist for all time}, for vertices $i\in V$ with $\varepsilon_i=1$, the fact that
$u_i(t)$ is upper bounded in $\mathbb{R}_{<0}$ is independent of the existence of hyperbolic discrete conformal structure
with combinatorial $\alpha$-curvature $\overline{R}$.
\end{remark}

\section{Open problems}\label{section 4}
\subsection{Rigidity of combinatorial $\alpha$-curvature}
The rigidity of combinatorial curvature for discrete conformal structures on polyhedral surfaces are usually studied under
the structure conditions (\ref{structure condition 1}) and (\ref{structure condition 2}) for the weight $\varepsilon: V\rightarrow \{0,1\}$,
which corresponds to the Euclidean and hyperbolic background geometry.
However, the rigidity of discrete conformal structures in the spherical background geometry is seldom studied, which corresponds to $\varepsilon=-1$.
Please refer to \cite{BBP, Luo2, MS} and others for some results.
It would be interesting to know a general version of the rigidity of discrete conformal structures on surfaces
in the spherical background geometry.

Another problem on rigidity is that the rigidity of combinatorial $\alpha$-curvature for different types of discrete conformal
structures on polyhedral surfaces are usually studied under the condition $\alpha \overline{R}\leq 0$,
where $\overline{R}$ is the prescribed combinatorial $\alpha$-curvature. Please refer to \cite{Ge-Xu 17, Ge-Xu 21, Xu I, Xu Z} and others.
By the proof of the rigidity, one can also generalize this condition to be a condition on the lower bound of the
eigenvalue of the $\alpha$-Laplace operator $\Delta_\alpha$. Please refer to \cite{Ge-Xu 21} (Section 5) or
the proof of Theorem \ref{Euclidean global rigidity theorem} in this paper.
A natural question on the rigidity of combinatorial $\alpha$-curvature is the following.

\begin{question}
What is the sharp lower bound of the eigenvalues of the $\alpha$-Laplace operator $\Delta_\alpha$ that ensures
the rigidity of the combinatorial $\alpha$-curvature? Does the sharp lower bound contain any geometric or topological
information on the triangulated surfaces?
\end{question}

\subsection{Prescribed combinatorial $\alpha$-curvature problem}
Prescribed Gaussian curvature on smooth surfaces is a fundamental problem in Riemannian geometry.
Prescribed combinatorial curvature on polyhedral surfaces also has lots
of applications in geometry, topology and applications, which asks the following question.
\\
\\
\textbf{Prescribed combinatorial curvature problem}
Suppose $M$ is a connected closed surface and  $V$ is a finite nonempty subset of $M$.
Which functions defined on $V$ are combinatorial curvatures of polyhedral metrics on $(M, V)$?
\\
\\
There are lots of important works on this problem for the classical combinatorial curvature
on polyhedral surfaces, which corresponds to $\alpha=0$ in this paper.
Thurston \cite{T1} gave a full characterization of the image of the classical curvature map for Thurston's circle packing
metrics on surfaces, which corresponds to $\varepsilon\equiv1$ and $\eta: E\rightarrow [0,1]$ in this paper.
Gu-Luo-Sun-Wu \cite{Gu1} and Gu-Guo-Luo-Sun-Wu \cite{Gu2} proved the fundamental
discrete uniformization theorem for Luo's vertex scaling by introducing a new definition of discrete conformality, which corresponds
to $\varepsilon\equiv 0$ in this paper.
In the case of $\alpha\neq 0$,
the authors \cite{Xu I,Xu Z,Xu Z 2} obtained some Kazdan-Warner type existence results for prescribed combinatorial $\alpha$-curvature problem in the setting of vertex scaling ($\varepsilon\equiv 0$) using the discrete conformal theory established in \cite{Gu1,Gu2}.
Some equivalent conditions for the solution of prescribed combinatorial $\alpha$-curvature problems on triangulated polyhedral surfaces
are given using the convergence of combinatorial $\alpha$-curvature flows
in \cite{GJ0, Ge-Jiang II, Ge-Jiang III, Ge-Jiang I, Ge-Xu 17, Ge-Xu 18, Ge-Xu 21} and others.

As the combinatorial $\alpha$-curvature can be taken as an approximation of the smooth Gaussian curvature on surfaces,
it would be interesting to get some existence results for the prescribed combinatorial $\alpha$-curvature problem
for generic $\alpha$, $\varepsilon$ and $\eta$.

\subsection{Combinatorial $\alpha$-curvature flows with surgery}
In Theorem \ref{main theorem 1}, we prove the convergence of combinatorial $\alpha$-Ricci flow by
extending the flow through removable singularities along the flow under the assumption that
there exists a discrete conformal structure with the prescribed combinatorial $\alpha$-curvature.
However, as noted in Remark \ref{Calabi flow not extendable}, the combinatorial $\alpha$-Calabi flow
can not be extended in this way.
On the other hand,
in the case of vertex scaling, Gu-Luo-Sun-Wu \cite{Gu1} and Gu-Guo-Luo-Sun-Wu \cite{Gu2}
introduced a different extension of Luo's combinatorial Yamabe flow \cite{Luo 1} via doing surgery along the
flow by edge flipping under the Delaunay condition. Using the surgery by edge flipping,
Gu-Luo-Sun-Wu \cite{Gu1} and Gu-Guo-Luo-Sun-Wu \cite{Gu2} proved the convergence of the combinatorial Yamabe flow
without the assumption of the existence of conformal factors with prescribed combinatorial curvature.
This method is then applied by Zhu-Xu \cite{Zhu Xu} to prove the convergence of the combinatorial Calabi flow with surgery.
The first author \cite{Xu I} and the authors \cite{Xu Z 2} further applied this surgery to combinatorial $\alpha$-curvature flows
and obtained some Kazdan-Warner type existence results for prescribed combinatorial $\alpha$-curvature problem in the setting of vertex scaling.
Motivated by these facts,
it is natural to introduce surgery by edge flipping under the weighted Delaunay condition along the combinatoiral $\alpha$-Ricci
flow and $\alpha$-Calabi flow for the discrete conformal structures in Definition \ref{edge length definition}.
We have the following conjecture on the combinatoiral $\alpha$-curvature flows with surgery.
\begin{conjecture}
Under the structure conditions (\ref{structure condition 1}) and (\ref{structure condition 2}) with $\varepsilon: V\rightarrow \{0, 1\}$,
the solution of  modified combinatoiral $\alpha$-Ricci flow and $\alpha$-Calabi flow with surgery
for discrete conformal structures in Definition \ref{edge length definition}
exists for all time and converges to the target combinatorial $\alpha$-curvature given by the prescribed combinatorial $\alpha$-curvature
problem.
\end{conjecture}

\bigskip

(Xu Xu) School of Mathematics and Statistics, Wuhan University, Wuhan 430072, P.R. China

E-mail: xuxu2@whu.edu.cn\\[2pt]

(Chao Zheng) School of Mathematics and Statistics, Wuhan University, Wuhan 430072, P.R. China

E-mail: 2019202010023@whu.edu.cn\\[2pt]


\begin{thebibliography}{50}
\setlength{\itemsep}{-2pt} \small

\bibitem{BPS} A. Bobenko, U. Pinkall, B. Springborn, \emph{Discrete conformal maps and ideal hyperbolic polyhedra}. Geom. Topol. 19 (2015), no. 4, 2155-2215.

\bibitem{BBP} Bowers, John C.; Bowers, Philip L.; Pratt, Kevin, \emph{Rigidity of circle polyhedra in the $2$-sphere and of hyperideal polyhedra in hyperbolic $3$-space}. Trans. Amer. Math. Soc. 371 (2019), no. 6, 4215-4249.

\bibitem{CMS} J. Cheeger, W. M\"{u}ller, R. Schrader, On the curvature of piecewise flat spaces, Comm. Math. Phys. 92 (3) (1984) 405-454.

\bibitem{Chow-Luo} B. Chow, F. Luo, \emph{Combinatorial Ricci flows on surfaces}, J. Differential Geom, Volume 63, No. 1 (2003), 97-129.

\bibitem{Dai-Ge} S. Dai, H. Ge, \emph{Discrete Yamabe flows with R-curvature revisited}. J. Math. Anal. Appl. 484 (2020), no. 1, 123681, 11 pp.

\bibitem{Ge1} H. Ge, \emph{Combinatorial methods and geometric equations}, Thesis (Ph.D.)-Peking University, Beijing. 2012. (In Chinese).

\bibitem{Ge2} H. Ge, \emph{Combinatorial Calabi flows on surfaces}, Trans. Amer. Math. Soc. 370 (2018), no. 2, 1377-1391.

\bibitem{Ge3} H. Ge, B. Hua, \emph{On combinatorial Calabi flow with hyperbolic circle patterns}, Adv. Math. 333 (2018), 523-538.

\bibitem{Ge-Hua 20} H. Ge, B. Hua, \emph{$3$-dimension combinatorial Yamabe flow in hyperbolic background geometry}, Trans. Amer. Math. Soc. 373 (2020),  no. 7, 5111-5140.

\bibitem{GJ0} H. Ge, W. Jiang, \emph{On the deformation of discrete conformal factors on surfaces}. Calc. Var. Partial Differential Equations  55  (2016),  no. 6, Art. 136, 14 pp.

\bibitem{Ge-Jiang II} H. Ge, W. Jiang, \emph{On the deformation of inversive distance circle packings, II}. J. Funct. Anal.  272 (2017),  no. 9, 3573-3595.

\bibitem{Ge-Jiang III} H. Ge, W. Jiang, \emph{On the deformation of inversive distance circle packings, III}. J. Funct. Anal.  272  (2017),  no. 9, 3596-3609.

\bibitem{Ge-Jiang I} H. Ge, W. Jiang, \emph{On the deformation of inversive distance circle packings, I}. Trans. Amer. Math. Soc. 327 (2019),  no. 9, 6231-6261.

\bibitem{Ge-Xu 14} H. Ge, X. Xu, \emph{Discrete quasi-Einstein metrics and combinatorial curvature flows in $3$-dimension},  Adv. Math. 267 (2014), 470-497.

\bibitem{Ge-Xu 16} H. Ge, X. Xu, \emph{$2$-dimensional combinatorial Calabi flow in hyperbolic background geometry}, Differential Geom. Appl. 47 (2016)  86-98.

\bibitem{Ge-Xu alpha-curvatures}  H. Ge, X. Xu, \emph{$\alpha$-curvatures and $\alpha$-flows on low dimensional triangulated manifolds}. Calc. Var. Partial Differential Equations 55 (2016), no. 1, Art. 12, 16 pp.

\bibitem{Ge-Xu 17} H. Ge, X. Xu, \emph{A discrete Ricci flow on surfaces with hyperbolic background geometry}, Int. Math. Res. Not. IMRN 2017, no. 11, 3510-3527.

\bibitem{Ge-Xu 18} H. Ge, X. Xu, \emph{On a combinatorial curvature for surfaces with inversive distance circle packing metrics}. J. Funct. Anal. 275 (2018), no. 3, 523-558.

\bibitem{Ge-Xu 21} H. Ge, X. Xu, \emph{A combinatorial Yamabe problem on two and three dimensional manifolds}, Calc. Var. Partial Differential Equations 60 (2021), no. 1, 20.

\bibitem{GXZ}  H. Ge, X. Xu, S. Zhang, \emph{Three-dimensional discrete curvature flows and discrete Einstein metrics}. Pacific J. Math. 287 (2017), no. 1, 49-70.

\bibitem{Glickenstein2011} D. Glickenstein, \emph{Discrete conformal variations and scalar curvature on piecewise flat two and three dimensional manifolds}, J. Differential Geom. 87 (2011), no. 2, 201-237.

\bibitem{Glickenstein} D. Glickenstein, \emph{Euclidean formulation of discrete uniformization of the disk}, Geom. Imaging Comput. 3 (2016), no. 3-4, 57-80.

\bibitem{Glickenstein-Thomas} D. Glickenstein, J. Thomas, \emph{Duality structures and discrete conformal variations of piecewise constant curvature surfaces}, Adv. Math. 320 (2017), 250-278.

\bibitem{Gu2} X. D. Gu, R. Guo, F. Luo, J. Sun, T. Wu, \emph{A discrete uniformization theorem for polyhedral surfaces II}, J. Differential Geom. 109 (2018), no. 3, 431-466.

\bibitem{Gu1} X. D. Gu, F. Luo, J. Sun, T. Wu, \emph{A discrete uniformization theorem for polyhedral surfaces}, J. Differential Geom. 109 (2018), no. 2, 223-256.

\bibitem{Guo} R. Guo, \emph{Local rigidity of inversive distance circle packing}, Trans. Amer. Math. Soc. 363 (2011) 4757-4776.

\bibitem{Kourimska thesis} H. Kou\v{r}imsk\'{a}, \emph{Polyhedral surfaces of constant curvature and discrete uniformization}. PhD thesis, Technische Universit\"{a}t Berlin, 2020.

\bibitem{Kourimska paper} H. Kou\v{r}imsk\'{a}, \emph{Discrete Yamabe problem for polyhedral surfaces}.
\href{https://arxiv.org/abs/2103.15693} {arXiv:2103.15693 [math.MG].}

\bibitem{Luo 1} F. Luo, \emph{Combinatorial Yamabe flows on surfaces}, Commun. Contemp. Math. 6 (2004), no. 5, 765-780.

\bibitem{Luo2} F. Luo, \emph{A characterization of spherical polyhedral surfaces}. J. Differential Geom. 74 (2006), no. 3, 407-424.

\bibitem{Luo3} F. Luo, \emph{Rigidity of polyhedral surfaces, III}, Geom. Topol. 15 (2011), 2299-2319.

\bibitem{MS}J. Ma, J. Schlenker, \emph{Non-rigidity of spherical inversive distance circle packings}, Discrete Comput. Geom. 47 (3) (2012) 610-617.

\bibitem{Pontryagin} L.S. Pontryagin, \emph{Ordinary differential equations}, Addison-Wesley Publishing Company Inc., Reading, 1962.

\bibitem{T1} W. Thurston, \emph{Geometry and topology of $3$-manifolds}, Princeton lecture notes 1976, \href{http://www.msri.org/publications/books/gt3m}{http://www.msri.org/publications/books/gt3m}.

\bibitem{WX} T. Wu, X. Xu, \emph{Fractional combinatorial Calabi flow on surfaces}, \href{https://arxiv.org/abs/2107.14102}{	arXiv:2107.14102 [math.GT].}

\bibitem{Xu 18} X. Xu, \emph{Rigidity of inversive distance circle packings revisited}, Adv. Math. 332 (2018), 476-509.

\bibitem{Xu 20} X. Xu, \emph{On the global rigidity of sphere packings on $3$-dimensional manifolds}, J. Differential Geom. 115 (2020), no. 1, 175-193.

\bibitem{Xu JFA 21} X. Xu, \emph{Combinatorial Calabi flow on $3$-manifolds with toroidal boundary}, J. Funct. Anal. 280 (2021), no. 11, 108990.

\bibitem{Xu A new proof} X. Xu, \emph{A new proof of Bowers-Stephenson conjecture}, Math. Res. Lett. 28 (2021), no. 4, 1283-1306.

\bibitem{Xu I} X. Xu, \emph{Parameterized discrete uniformization theorems and curvature flows for polyhedral surfaces, I}. \href{https://arxiv.org/abs/1806.04516} {arXiv:1806.04516v2 [math.GT].}

\bibitem{Xu 21} X. Xu, \emph{Rigidity and deformation of discrete conformal structures on polyhedral surfaces}. \href{https://arxiv.org/abs/2103.05272} {arXiv:2103.05272 [math.DG].}

\bibitem{Xu Z} X. Xu, C. Zheng, \emph{Parameterized discrete uniformization theorems and curvature flows for polyhedral surfaces, II}. \href{https://arxiv.org/abs/2103.16077} {arXiv:2103.16077 [math.GT].} accepted by Trans. Amer. Math. Soc. 2021.

\bibitem{Xu Z 2} X. Xu, C. Zheng, \emph{Prescribing discrete Gaussian curvature on polyhedral surfaces}. \href{https://arxiv.org/abs/2110.12326} {arXiv:2110.12326 [math.DG].} accepted by Calc. Var. Partial Differential Equations. 2022.

\bibitem{Zeng-Gu} W. Zeng, X. Gu, \emph{Ricci flow for shape analysis and surface registration}, Springer Briefs in Mathematics. Springer, New York (2013).

\bibitem{ZGZLYG} M. Zhang, R. Guo, W. Zeng, F. Luo, S.T. Yau, X. Gu, \emph{The unified discrete surface Ricci flow}, Graphical Models 76(2014), 321-339.

\bibitem{Zhu Xu} X. Zhu, X. Xu, \emph{Combinatorial Calabi flow with surgery on surfaces}, Calc. Var. Partial Differential Equations 58 (2019), no. 6, Paper No. 195, 20 pp.


\end{thebibliography}
\end{document}